\documentclass[11pt,a4paper]{article}
\usepackage{dingbat}
\usepackage{amsmath,amssymb,latexsym}
\usepackage{theorem}
\usepackage{overpic}

\newtheorem{theorem}{Theorem}
\newtheorem{lemma}{Lemma}

\newtheorem{corollary}{Corollary}
\newtheorem{proposition}{Proposition}
{\theorembodyfont{\rmfamily} 

\newtheorem{remark}{Remark}}
{\theorembodyfont{\slshape} }

\newcommand{\res}{\mathop{\rm res}}

\newcommand{\field}[1]{\mathbb{#1}}
\newcommand{\A}{\field{A}}
\newcommand{\D}{\field{D}}
\newcommand{\R}{\field{R}}
\newcommand{\Z}{\field{Z}}
\newcommand{\N}{\field{N}}
\newcommand{\C}{\field{C}}
\newcommand{\Q}{\field{Q}}

\newcommand{\T}{\field{T}}
\newcommand{\BB}{{\mathcal B}}
\newcommand{\CC}{{\mathcal C}}

\newcommand{\FF}{{\mathcal S}}

\renewcommand{\AA}{{\mathcal A}}
\newcommand{\const}{{\rm const}}

\renewcommand{\Re}{\mathop{\rm Re}}
\renewcommand{\Im}{\mathop{\rm Im}}

\newcommand{\isdef}{\stackrel{\text{\tiny def}}{=}}

\DeclareRobustCommand{\qed}{%
\ifmmode % if math mode, assume display: omit penalty etc.
\else \leavevmode\unskip\penalty9999 \hbox{}\nobreak\hfill \fi
\quad\hbox{\qedsymbol}}
\newcommand{\openbox}{\leavevmode
\hbox to.77778em{%
\hfil\vrule
\vbox to.675em{\hrule width.6em\vfil\hrule}%
\vrule\hfil}}
\newcommand{\qedsymbol}{\openbox}
\newcommand{\proofname}{Proof}
\newenvironment{proof}[1][\proofname]{\par
\normalfont \trivlist \item[\hskip\labelsep   \itshape #1. ]
\ignorespaces
}{%
\qed \endtrivlist } %%%%%%%%%%%%%%%%%%%

\def\XXint#1#2#3{{\setbox0=\hbox{$#1{#2#3}{\int}$}
\vcenter{\hbox{$#2#3$}}\kern-.5\wd0}}

\title{Asymptotics of orthogonal polynomials with respect to an
analytic weight with algebraic singularities on the circle
}%

\author{A.\ Mart\'{\i}nez-Finkelshtein\footnote{Corresponding author.
E-mail: \texttt{andrei@ual.es}}\\ University of Almer\'{\i}a, Spain \and
K. T.-R. McLaughlin\\ University of Arizona, Tuckson, USA \and
E.\ B.\ Saff \\Vanderbilt University, USA }%

\date{}%

\begin{document}
\maketitle

\begin{abstract}
Strong asymptotics of polynomials orthogonal on the unit circle with
respect to a weight of the form
$$
W(z) =  w(z)\, \prod_{k=1}^m |z-a_k|^{2\beta _k} \,, \quad |z|=1 \,,
\quad |a_k|=1, \quad \beta _k>-1/2, \quad k=1, \dots, m,
$$
where
%$|a_k|=1 $ and $\beta _k>-1/2$, $k=1, \dots, m$, and where
$w(z)>0$ for $|z|=1$ and can be extended as a holomorphic and
non-vanishing function to an annulus containing the unit circle. The
formulas obtained are valid uniformly in the whole complex plane. As
a consequence, we obtain some results about the distribution of
zeros of these polynomials, the behavior of their leading and
Verblunsky coefficients, as well as give an alternative proof of the
Fisher-Hartwig conjecture about the asymptotics of Toeplitz
determinants for such type of weights. The main technique is the
steepest descent analysis of Deift and Zhou, based on the matrix
Riemann-Hilbert characterization proposed by Fokas, Its and Kitaev.
\end{abstract}

\maketitle

% ----------------------------------------------------------------

\section{Introduction and statement of the main results} \label{section:intro}

Let us first set some notation that will be widely used in what
follows. We denote by $\T $ the unit circle on the complex plane
$\C$ (circle of radius $1$ centered at the origin), and $\D \isdef
\{z\in \C:\, |z|<1\}$ is the open unit disc. If $r<1$, let
$\A_r\isdef \{z\in \C:\, r<|z|<1/r \}$. Every oriented Jordan curve
or arc $\gamma$ induces naturally the notions of left (or ``$+$'')
and right (or ``$-$'') sides of $\gamma$. We also denote by a
subindex ``$+$'' (respectively, ``$-$'') the left (respectively,
right) boundary values of functions on $\gamma $.

An integrable non-negative function $W$ defined on $\T$ is called a
\emph{weight} if
\begin{equation}\label{szegoCOndition}
\int_{\T} W(z)\, |dz| >0\,.
\end{equation}
For each weight $W$ there exists a unique sequence of polynomials
$\varphi_n$ (called \emph{Szeg\H{o} polynomials}), orthonormal with
respect to $W$, satisfying $\varphi_n (z)=\kappa_n z^n+\text{lower
degree terms}$, $\kappa_n >0$, and
\begin{equation}\label{orthogonalityConditions}
\oint_{\T} \varphi_n(z)  \overline{\varphi_m(z)}\, W(z)
|dz|=\delta_{mn}\,.
\end{equation}
We denote by $\Phi_n(z) \isdef  \varphi_n(z)/\kappa _n$ the
corresponding monic orthogonal polynomials. It is well known that
they satisfy the Szeg\H{o} recurrence
\begin{equation}\label{recurrence}
\Phi_{n+1}(z)=z \Phi_n(z)-\overline{\alpha _n}\, \Phi_n^*(z)\,,
\quad \Phi_0(z)\equiv 1\,,
\end{equation}
where we use the standard notation $\Phi_n^*(z)\isdef z^n
\overline{\Phi_n(1/\overline{z})}$. The parameters $\alpha
_n=-\overline{\Phi_{n+1}(0)}$ are called \emph{Verblunsky
coefficients} (also \emph{reflection coefficients} or \emph{Schur
parameters}) and satisfy $\alpha _n \in \D$ for $n=0, 1, 2, \dots$
(see \cite{Simon05b} for details).

If the weight $W$ satisfies the condition
\begin{equation}\label{condition_Szego}
\int_{\T} \log W(z) \, |dz|>-\infty\,,
\end{equation}
then the \emph{Szeg\H{o} function} of $W$ (see e.g.\ \cite[Ch.\ X,
\S 10.2]{szego:1975}),
\begin{equation}\label{standardSzego}
D (W; z) \isdef \exp\left( \frac{1}{4\pi }\,\int_0^{2\pi} \log
W(e^{i \theta}) \, \frac{e^{i \theta}+z}{e^{i \theta}-z}\, d\theta
\right)\,,
\end{equation}
can be defined. This function is piecewise analytic and
non-vanishing, defined for $ z\notin \T$, and we will denote by
$D_{\rm i}$ and $D_{\rm e}$ its values for $|z|<1$ and $|z|>1$,
respectively. It is easy to verify that
\begin{equation}\label{symmetry}
\overline{D_{\rm i}\left(W; \frac{1}{\overline{z}}
\right)}=\frac{1}{D_{\rm e}(W; z)}\,, \quad |z|>1\,,
\end{equation}
and in particular, $W(z)=|D_{\rm e}(W; z)|^{-2}$ for $z\in \T $. The
role of the Szeg\H{o} function in the description of the asymptotic
behavior of the orthogonal polynomials is well known; for instance,
$$
\lim_{n\to \infty} \frac{\Phi_n(z)}{z^n}=   \frac{D_{\rm
e}(W;z)}{D_{\rm e}(W;\infty)}
$$
uniformly on each compact set in the exterior of the unit disk.

In this paper we focus on weights of a specific form. Assume that
$w$ is a strictly positive function defined on the unit circle,
which can be extended as a holomorphic and non-vanishing function to
an annulus $\A_\rho$ ($0<\rho<1$); in fact, a ``canonical'' form for
such a function is $w(z)=|f(z)|^2$, $z\in \T$, where $f$ is
holomorphic and non-vanishing in $\A_\rho$. For points $a_k \in \T $
and values $\beta _k>-1/2$, $k=1, \dots, m$, we define the following
weight:
\begin{equation}\label{THE_weight}
W(z) \isdef  w(z)\, \prod_{k=1}^m |z-a_k|^{2\beta _k} \,, \quad z\in
\T \,.
\end{equation}
It can have zeros or blow up at $a_k$'s, but still conditions
\eqref{szegoCOndition} and \eqref{condition_Szego} are satisfied. We
are interested in the asymptotic behavior of the corresponding
sequences $\{\Phi_n (z)\}$, $\{ \kappa_n\}$ and $\{ \alpha _n\}$
when $n\to \infty$ and $z \in \C$.

The case when all $\beta _k=0$ (that is, when $W$ is a positive
analytic weight on $\T $) has been studied in
\cite{math.CA/0502300}, where a canonical representation of the
corresponding orthogonal polynomials in terms of iterates of the
Cauchy transform of the \emph{scattering function} of $W$,
\begin{equation}\label{def_Scattering}
    \FF(W;z) \isdef D_{\rm i}(W; z)D_{\rm e}(W; z)\,,
\end{equation}
was derived. Also, detailed asymptotic formulas were obtained. A
characterizing feature of this case is that the zeros of $\Phi_n$'s
stay away from $\T$, clustering (with a possible exception of a
$o(n)$ number of them) at an inner circle determined by the analytic
continuation of $D_{\rm e}$.

When any $\beta _k\neq 0$,  $W$ is no longer a positive and analytic
weight on $\T $, and in this situation the majority of the zeros of
the Szeg\H{o} polynomials are attracted by the unit circle. One of
the main goals of the paper is to provide asymptotic formulas for $
\Phi_n $'s valid uniformly on the whole complex plane, and in
particular, in a neighborhood of $a_k$'s. For partial results in the
case when all $\beta_k=1$ see \cite{Bello06}.

In order to state the main results we need to introduce a new piece
of notation.

Let  $\AA \isdef \{a_1, \dots, a_m \}$; for the sake of brevity
hereafter we call these points generically as ``singularities'' of
the weight, although $W$ is regular at $a_k$'s for integer $\beta
_k$'s. We fix for what follows $0 <\delta <1-\rho $, such that
additionally $\delta <\frac{1}{3}\min_{i\neq j} |a_i-a_j|$, and
denote
\begin{equation}\label{notationBC}
    \BB_k\isdef \{z\in \C:\, |z-a_k|<\delta \}\,, \quad  \CC_k\isdef \{z\in \C:\,
    |z-a_k|=\delta\}
    \,, \quad k=1, \dots, m\,,
\end{equation}
as well as $B\isdef \cup_{k=1}^m \BB_k$. Furthermore, given a subset
$X\subset \C$ and a value $a\in \C$ we will use the standard
notation $ a \cdot X =\{a x:\, x\in X \}$; consistently, $\AA \cdot
X\isdef \cup_{k=1}^m (a_k\cdot X)$.

Under assumptions on $W$ in \eqref{THE_weight}, both $D_{\rm
i}(W;z)$ and $D_{\rm e}(W;z)$ admit an analytic extension across $\T
\setminus  \AA$; we keep the same notation for these analytic
continuations. More precisely, $D_{\rm i}(W; z)$ is holomorphic in
$\{ z\in\C:\, |z|<1/ \rho \} \setminus \left(\AA \cdot [1,1/\rho )
\right)$, and $D_{\rm e}(W; z)$ is holomorphic in $\{ z\in\C:\,
|z|>\rho \} \setminus \left(\AA \cdot (\rho,1] \right)$; see Section
\ref{sec:scattering} for a detailed discussion.

First we describe the asymptotic behavior of $\Phi_n$'s away from
the singular points of the weight:
\begin{theorem}
\label{thm1} For monic orthogonal polynomials $\Phi_n$ corresponding
to the weight $W$ given in \eqref{THE_weight} there exist a complete
asymptotic expansion
\begin{equation}\label{completeExpansionofPhi}
\Phi_n(z)= \sum_{k=0}^\infty \frac{\mathfrak f_k(z)}{n^k}
\end{equation}
valid uniformly in $\C$. Each $\mathfrak f_k(z)$ is a piecewise
analytic function, holomorphic in each domain specified below. In
particular, there exist constants $\vartheta_k \in \T$, $k=1, \dots,
m$, defined by formula \eqref{def_new_constant} below, such that:
\begin{enumerate}
\item[(i)] formula
\begin{equation*}%\label{asymptForPolynsInSigma0}
\Phi_n(z)=  \frac{ D_{\rm i}(W; 0)}{  D_{\rm i}(W; z)}\,
\frac{1}{n}\, \left( \sum_{k=1}^m \dfrac{\beta _k \vartheta_k
 }{ a_k-z }\, a_k^{n+1}    + O\left(\frac{1}{n }\right) \right)
\end{equation*}
holds uniformly on every compact subset of $\D $;
\item[(ii)] formula
\begin{equation*}%\label{asympExprinOmegaPlus}
 \begin{split}
\Phi_n(z) = &  z^n \, \frac{D_{\rm e}(W;z)}{D_{\rm e}(W;\infty)}\,
\left(1+ \frac{1}{n}\, \sum_{k=1}^m \frac{a_k \beta _k^2  }{ a_k-z }
+\mathcal O \left(\frac{1}{n^2}\right)\right) \\ & + \frac{D_{\rm
i}(W;0)}{ D_{\rm i}(W;z)}\, \left(\frac{1}{n}\, \sum_{k=1}^m \frac{
\beta _k \vartheta_k}{ a_k- z}\,  a_k^{n+1}  +\mathcal O
\left(\frac{1}{n^2}\right) \right)
 \end{split}
\end{equation*}
holds uniformly on every compact subset of  $\A_\rho \setminus
\left(B \cup (\AA \cdot (\rho ,1/\rho )\right) $;
\item[(iii)] formula
\begin{equation*}%\label{asympExprinOmegaInfty}
 \begin{split}
\Phi_n(z)& =z^n \, \frac{D_{\rm e}(W;z)}{D_{\rm e}(W;\infty)}\,
\left(1+ \frac{1}{n}\, \sum_{k=1}^m \frac{a_k \beta _k^2  }{  a_k-z
} +\mathcal O \left(\frac{1}{n^2}\right)\right) \,.
 \end{split}
\end{equation*}
holds uniformly on every compact subset of $\C\setminus \overline{\D
}$.
\end{enumerate}

\end{theorem}
\begin{remark}
It is well known (see \cite{MR2001g:42050} as well as Section
\ref{sec:RH} below) that further terms $\mathfrak  f_k$ of the
expansion \eqref{completeExpansionofPhi} can be obtained by nested
contour integration and the calculus of residues. However, the
difficulty of the computation increases with $k$, and we limit our
attention to the leading nontrivial terms of
\eqref{completeExpansionofPhi}.
\end{remark}

\begin{corollary} \label{cor:zeros} %\doubt
For every compact set $K\subset \D$ there exists $N=N(K)\in \N$ such
that for every $n\geq N$, every $\Phi_n$ has at most $m-1$ zeros on
$K$.
\end{corollary}

However, the global behavior of these ``spurious'' $m-1$ zeros can
be complicated; we describe their limiting set below (Theorem
\ref{thm:szabados}).
%%%%%%%%%%%%%%%%%%%%%%%%%%%%%%%%%%%%%%%%%%%%%%%%%%%%%%%%%%
\begin{figure}[htb]
\centering \begin{overpic}[scale=1.1]{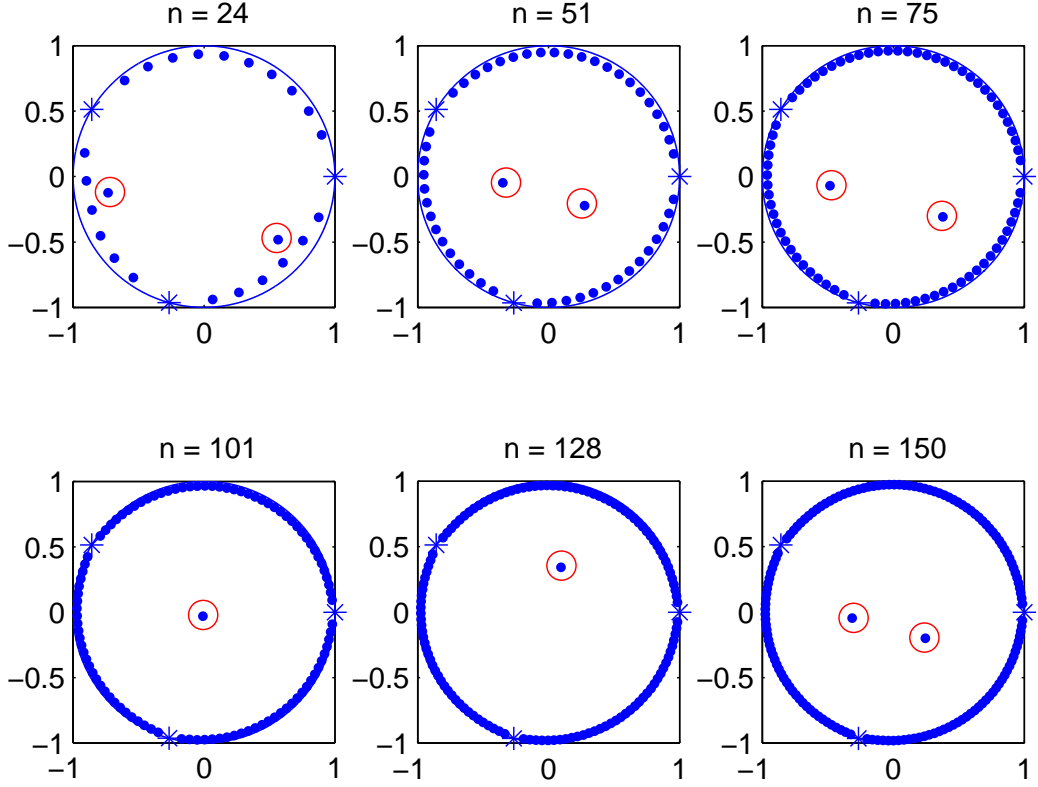}
\end{overpic}
\caption{Zeros of $\Phi_n$ for several values of $n$ with
$W(z)=|(z-1) (z-a)(z-a^2)|^4$, $z\in \T$, where $a=\exp(\pi i
\sqrt{2})$. Centers of the circles are the zeros of $ \sum_{k=1}^m
 \beta _k \vartheta_k a_k^{n+1}/(
  a_k-z) $.}\label{fig:spurious}
\end{figure}
%%%%%%%%%%%%%%%%%%%%%%%%%%%%%%%%%%%

In order to formulate the asymptotic behavior of the Szeg\H{o}
polynomials in a neighborhood of each singular point $a_k$ we define
an auxiliary function: if $\zeta ^{1/2}$ stands for the main branch
of the square root in $ \C\setminus (-\infty, 0]$ (that is, $\zeta
^{1/2}>0$ for $\zeta >0$), and $J_\nu $ are the Bessel functions of
the first kind, then for $\beta >-1/2$ set
\begin{equation}\label{def_H_for_local}
    \mathcal H(\beta ; \zeta )\isdef \begin{cases}
    e^{-2\pi i \beta } \zeta ^{1/2} \left( i J_{\beta
    +1/2}(\zeta )+ J_{\beta -1/2}(\zeta )\right)\,, & \text{if $\zeta$ is in the second quadrant,}
    \\ \zeta ^{1/2} \left( i J_{\beta
    +1/2}(\zeta )+ J_{\beta -1/2}(\zeta )\right)\,, &
    \text{otherwise.}
    \end{cases}
\end{equation}

With the notations introduced in \eqref{notationBC} we have the
following result about local behavior of the polynomials at the
singularities of the weight:
\begin{theorem}
\label{thm:local_behavior} Let $k\in \{1, \dots, m\}$, $|z-a_k|\leq
\delta $. For $z\in \BB_k$ define
\begin{equation}\label{def_zeta_local}
\zeta_n(z) \isdef -i \,\frac{n}{2}\, \log\left(\frac{z}{a_k}\right)
\,,
\end{equation}
where we take the principal branch of the logarithm. Then for $z\in
\BB_k$,
\begin{equation}\label{Local_asymptotics}
\begin{split}   \Phi_n(z) &= \sqrt{\frac{\pi   } { 2}} \,e^{\pi i \beta
_k/2}\,   \dfrac{D_{\rm e}(W;z)}{D_{\rm e}(W;\infty)}\,  \left(
a_k\, z\right)^{n/2} \, \mathcal H\left(\beta_k ; \zeta_n(z)
\right)\, \left(
1+\mathcal O\left(\frac{1}{n}\right)\right)\,, %\quad z\in \BB_k\,,
\end{split}
\end{equation}
where we take the principal branch of the square root in the
described neighborhood of $z=a_k$. The $ \mathcal O\left(1/n\right)$
term in \eqref{Local_asymptotics} is uniform in the closed disk
$\overline{\BB_k}$.
\end{theorem}
\begin{remark} \label{remark:local}
It will be shown in Section \ref{sec:asymptotics} that  $D_{\rm
e}(W;z) \mathcal H(\beta_k ; \zeta_n(z) )$ is a holomorphic function
in a neighborhood of $z=a_k$.
\end{remark}
\begin{remark}
In the particular case of $\beta _k=0$ (a removable singularity) we
have
$$
\mathcal H(0 ; \zeta ) = \zeta^{1/2} \left( i   J_{ \frac{1}{2}}
(\zeta) +
 J_{  -\frac{1}{2}} (\zeta) \right) = \sqrt{\frac{\pi   } { 2}} \, e^{i\zeta }\,,
$$
and \eqref{Local_asymptotics} takes the form
\begin{equation*}
\begin{split}   \Phi_n(z) &=  \frac{z^n D_{\rm e}(W;z)}{D_{\rm e}(W;\infty) }\, \, \left( 1+\mathcal
O\left(\frac{1}{n}\right)\right)\,,
\end{split}
\end{equation*}
(cf.\ Theorem \ref{thm1}, \emph{(ii)--(iii)}).
\end{remark}
\begin{remark}
The steepest descent analysis used for the proof of the theorems
above allows to find further terms of the asymptotic expansion for
the polynomials. However, for the sake of simplicity we decided to
restrict our attention to the leading terms bearing already
non-trivial information about the main parameters and the zeros, as
it will be shown next.
\end{remark}

For a weight  $W$ on $\T$ let us define the geometric mean
\begin{equation}\label{def_G}
G[W]\isdef   \exp\left(\frac{1}{2\pi} \int_{0}^{2\pi}
 \log\left( W\left(e^{i\theta} \right)\, d\theta \right)
 \right)\,.
\end{equation}
If $W$ is given by \eqref{THE_weight}, then $G[W]=G[w]= D_{\rm
i}^2(W;0)$ (see formulas \eqref{tau}--\eqref{tauAlternativeBis}
below).

With respect to the Verblunsky and leading coefficients of the
orthonormal polynomials we have
\begin{theorem}
\label{thm:coefficients} For weight $W$ given in \eqref{THE_weight}
there exist constants $\vartheta_k \in \T$, $k=1, \dots, m$,
introduced in \eqref{def_new_constant} below, such that the
Verblunsky coefficients $\alpha _n$ defined by the recurrence
\eqref{recurrence} satisfy
\begin{equation}\label{result_for_Verblunsky}
\alpha _n =  -\frac{1}{n}\, \sum_{k=1}^m \frac{\beta _k}{
\vartheta_k }\, \frac{1}{a_k^{n+1}}  + \mathcal O\left(\frac{1}{n^2
} \right) \,, \quad n \to \infty\,.
\end{equation}
For the leading coefficients $\kappa_n$ of the orthonormal
polynomials $\varphi_n$ the following asymptotic formula holds:
\begin{equation}\label{result_for_kappa}
\kappa _{n-1}^2= \frac{1}{G[2\pi w]}\,\left( 1 - \frac{1}{n}\,
\sum_{k=1}^m \beta _k^2+\mathcal O\left( \frac{1}{n^2}\right)
\right) \,, \quad n\to \infty.
\end{equation}
\end{theorem}
\begin{remark}
Obviously formula \eqref{result_for_kappa} makes sense for $n\geq
\sum_{k=1}^m \beta _k^2$, and exhibits, at least asymptotically, the
growing character of the leading coefficients $\kappa _n$.
\end{remark}

A direct consequence of \eqref{result_for_kappa} is the asymptotic
behavior of the Toeplitz determinants related to the weight $W$. If
we define the moments
$$
d_k\isdef \oint_{z\in \T } z^{-n} W(z) |dz|\,,
$$
then the Toeplitz determinants are
\begin{equation}\label{def_Toeplitz_dets}
\mathcal D_n(W)\isdef \det \left[ \left(
d_{j-i}\right)_{i,j=0}^n\right]\,.
\end{equation}
\begin{theorem} \label{thm:Fisher}
Under the assumption above there exists a constant $\varkappa$
depending on $W$ such that
\begin{equation}\label{asympt_Tplitz}
\begin{split}
\mathcal D_n(W) & =\varkappa\,  \left( G[2\pi w] \right)^n\,
  n^{\sum_{k=1}^m \beta _k^2} \, \left( 1+o (1)\right) \,, \quad n \to \infty\,.
\end{split}
\end{equation}
\end{theorem}
This formula is in accordance with the well known Fisher-Hartwig
conjecture (see e.g.\ \cite{Basor91}).
\begin{remark}
As it was mentioned above, in the case considered, $G[  w]=G[  W]$,
so we may replace $w$ by $W$ in the right hand side of both formulas
\eqref{result_for_kappa} and \eqref{asympt_Tplitz}.
\end{remark}

Let us discuss now how Theorems \ref{thm1}--\ref{thm:local_behavior}
reveal the behavior of the zeros of the polynomials $\Phi_n$.
Qualitatively we can describe the picture as follows: the vast
majority of the zeros will approach the circle $\T $ regularly and
radially uniformly along level curves $\Gamma_n$ defined below, that
are asymptotically close to circles centered at the origin with
radius $n^{-1/n}$. Singularities of the weight  exert the following
influence on the zeros of $\Phi_n$: those closest to points $a_k$
converge to $\T$ faster than the rest (their absolute value behaves
like $c^{-1/n}$, with $c>1$ depending on $\beta _k$), and either
leaving a ``gap'' around $a_k$ (if $\beta _k>0$) or approaching it
radially (if $\beta _k<0$). Furthermore, a bounded number of
``spurious'' zeros may wander inside the unit disk; each $\Phi_n$
will have at most $m-1$ of these zeros (Corollary \ref{cor:zeros}),
and  their global behavior will depend in particular on the relative
positions of points $a_k$.

Let
\begin{equation}\label{def_LevelCurve}
\Gamma_n\isdef \left\{ z\in \C:\,  |z|^n |\FF(W;z)| = \frac{1}{n}\,
\left| \sum_{k=1}^m \frac{\beta _k \vartheta_k a_k ^{n+1}}{z- a_k }
\right|
 \right\}\,,
\end{equation}
where $\mathcal S(W;z)$ is the scattering function defined in
\eqref{def_Scattering}. These curves are well defined: although
$\mathcal S(W;z)$ is multivalued in a neighborhood of $\T$, by
\eqref{boundaryValueF_onCuts} below, $|\mathcal S(W;z)|$ is
single-valued and positive in $\A_\rho $. Furthermore, given an
analytic function $f$ and $\varepsilon>0$ denote
\begin{equation}\label{nullSet}
    \mathcal Z(f) \isdef \{z:\, f(z)=0 \}\,,
    \quad \mathcal Z_\varepsilon (f )
    \isdef  \left\{ z\in \C:\, \min_{t\in \mathcal Z (f )}
|t-z|<\varepsilon \right\} \,,
\end{equation}
and let
$$
\Gamma_n(\varepsilon )\isdef \left\{z\in \C:\, \min_{t\in \Gamma_n}
|t-z|<\varepsilon \right\}\,,
$$
be the $\varepsilon$-neighborhood of $\Gamma_n$.

Let also
\begin{equation}\label{defTheta}
\theta_1=1, \quad \text{and} \quad
\theta_k=\frac{1}{2\pi}\,\left(\arg a_k - \arg a_1 \right),
  \qquad  k=2, \dots, m\, ,
\end{equation}
and
\begin{equation}\label{G_n}
   \mathcal R_n(z)\isdef \sum_{k=1 }^m\frac{\beta _k \vartheta_k  }{z- a_k }\, e^{2\pi i (n+1)
\theta_k}\,.
\end{equation}

\begin{figure}[htb]
\centering \begin{overpic}[scale=1.1]{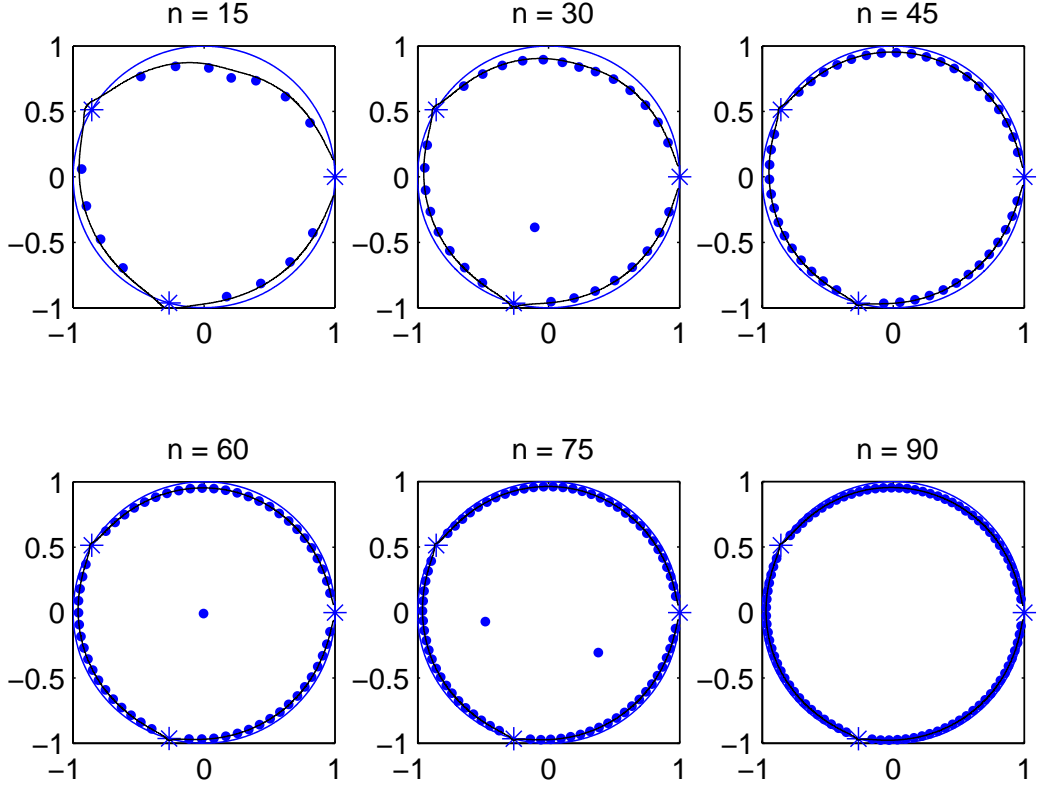}
\end{overpic}
\caption{Zeros of $\Phi_n$ for $n=15, 30, \dots, 90$  with
$W(z)=|(z-1) (z-a)(z-a^2)|^4$, $z\in \T$, where $a=\exp(\pi i
\sqrt{2})$. Points $1$, $a$ and $a^2$ are indicated with asterisks.
For comparison with the prediction of Theorem \ref{thm:clock} we
plot in each case the level curve $\Gamma_n$ defined in
\eqref{def_LevelCurve}.}\label{fig:3poles4_ratcurves}
\end{figure}

\begin{theorem}
\label{thm:clock} There exists $0<\varepsilon <\delta$ such that for
all sufficiently large $n\in\N$ every ``pie slice'' of the form
\begin{equation*}
\begin{split}
 \left\{z\in \Gamma_n(\varepsilon):\,  \; \frac{\alpha +2 k_1 \pi}{n}<\arg\left(z
\right)<\frac{\alpha  +2 k_2 \pi}{n} \right\}\subset \A_\rho
\setminus (B \cup Z_\varepsilon (\mathcal R_n ))\,,
\end{split}
\end{equation*}
with appropriately chosen $\alpha \in \R$, contains exactly
$k=k_2-k_1 $ zeros of $\Phi_n$, $z^{(n)}_1,\dots, z^{(n)}_k$,
satisfying
\begin{equation}\label{asymptModEquidistr}
|z^{(n)}_i|= 1-\frac{\log(n)}{n}+\mathcal O\left(\frac{1}{n}\right)
\,,
\end{equation}
and
\begin{equation}\label{asymptArgEquidistr}
\arg\left(z^{(n)}_{i+j}\right)-\arg\left(z^{(n)}_i\right)=\frac{2\pi
j}{n}+\mathcal O \left( \frac{1}{n^2} \right)\,.
\end{equation}

\end{theorem}
See Figure \ref{fig:3poles4_ratcurves} for an illustration of the
statements of Theorem \ref{thm:clock}.

%%%%%%%%%%%%%%%%%%%%%%%%%%%%%%%%%%%%%%%%%%%%%%%%%%%%%%%%%%
%%%%%%%%%%%%%%%%%%%%%%%%%%%%%%%%%%%

We can also be more precise about the accumulation set of zeros of
$\Phi_n$'s,
\begin{equation}\label{def:accumulationSet}
    Z \isdef \bigcap_{k\geq 1} \overline{\bigcup_{n\geq k} \mathcal Z( \Phi_n)
    }\,,
\end{equation}
on compact sets of $\D $. Assertion \emph{(i)} of Theorem \ref{thm1}
shows that the structure of $Z\cap \D  $ depends on the relative
positions of $a_1, \dots, a_m$ on $\T  $. Without loss of generality
we may assume that $\theta_1=1, \theta_2, \dots, \theta_v$ is the
maximal subset of $\{ \theta_1, \dots, \theta_m\}$ linearly
independent over the rational numbers $\Q$. Then there exist unique
$r_{kj}\in \Q$, $k=1, \dots, m$, $j=1, \dots, v$, such that
$$
\theta_k=\sum_{j=1}^v r_{kj}\, \theta_j\,, \quad \ k=1, \dots, m\,,
$$
(obviously, $r_{kj}=\delta_{kj}$, for $1\leq k \leq v$). Let
$q_1=1$, and if $m\geq 2$, $ 0\leq p_k< q_k\in \Z$, $ k=2, \dots, m$
be such that
$$
r_{k1} \equiv \frac{p_k}{q_k} \mod \Z, \quad   k=2, \dots, m\,.
$$

\begin{theorem}
\label{thm:szabados} Let $t\in Z\cap \D  $. If $v=1$ (that is, if
all $\theta_k\in \Q$), then there exist $0\leq s_k<q_k$, $s_k\in
\Z$, $k=1,\dots,m$, such that
\begin{equation}\label{description_Szabados}
\sum_{k=1}^m \dfrac{\beta _k \vartheta_k
 }{ a_k-t } \, e^{2\pi i  s_k/q_{k} }=0\,.
\end{equation}

If $v\geq 2$, then additionally there exist $X_2, \dots, X_v\in \R$
such that
\begin{equation}\label{description_Szabados2}
\sum_{k=1}^m \dfrac{\beta _k \vartheta_k
 }{ a_k-t } \, e^{2\pi i \left( s_k/q_{k}+\sum_{j=2}^v r_{kj} \, X
_j\right)}=0\,.
\end{equation}
Moreover, every point in $\D$ that belongs to the manifold given by
solutions of \eqref{description_Szabados2} when $X_2, \dots, X_v$
vary in $\R$, is an accumulation point of the zeros of $\{
\Phi_n\}$.
\end{theorem}

\begin{corollary}\label{cor:cluster}
\begin{enumerate}
\item[(i)] If $v=1$, then $Z\cap \D$ is a discrete set of
a finite number of points.

\item[(ii)] If $v=2$, then $Z\cap \D$ is an algebraic curve of degree $\leq m$.
In particular, if $v=m=2$, then $Z\cap \D $ is either a circular arc
(if $|\beta _1|\neq |\beta _2|$), or a diameter in $\D $ formed by
the perpendicular bisector of the segment joining $a_1$ with $a_2$
(if $|\beta _1|=|\beta _2|$).

\item[(iii)] If $v> 2$, then $Z\cap \D$ is a two-dimensional domain
bounded by algebraic curves.
\end{enumerate}

\end{corollary}

\begin{remark}
From the method of proof it follows in fact that the zeros of the
orthogonal polynomials will be equidistributed on the manifolds
described above.
\end{remark}
%%%%%%%%%%%%%%%%%%%%%%%%%%%%%%%%%%%%%%%%%%%%%%%%%%%%%%%%%%
\begin{figure}[htb]
\centering \begin{overpic}[scale=0.85]{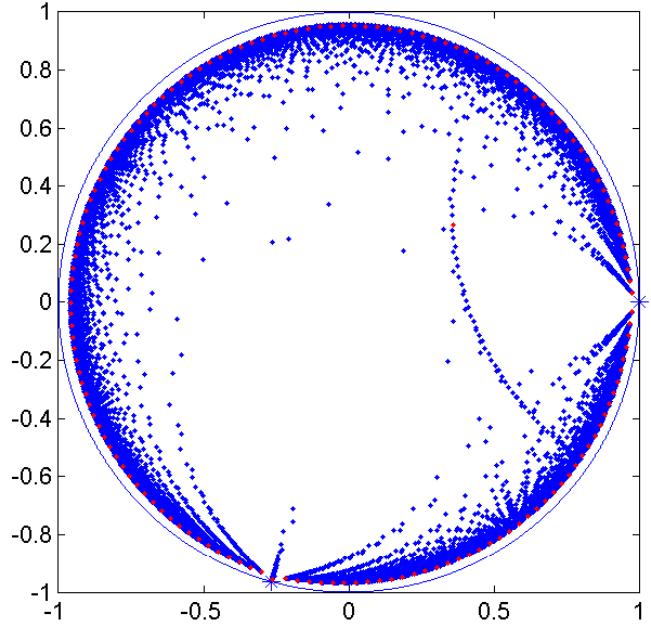}
\end{overpic}
\caption{Zeros of $\Phi_n$ for $n=1, 2, \dots, 150$  with
$W(z)=|z-1|^{1/3} |z-a|^{-2/3}$, $z\in \T$, where $a=\exp(\pi i
\sqrt{2})$. Points $1$ and $a$ are indicated with asterisks. Observe
different features of the behavior of the zeros described in the
text: the zero of the weight at $z=1$ repels the zeros of the
polynomial, while the singularity at $z=a$ attracts a zero of
$\Phi_n(z)$ that approaches $z=a$ radially (cf.\ Remark
\ref{remark:localBehavior}). Also the arc of the circle along which
the zeros of $\Phi_n$ remaining inside $\D$ cluster is clearly
visible (cf.\ Corollary \ref{cor:cluster}).}\label{fig:2poles1irrat}
\end{figure}
%%%%%%%%%%%%%%%%%%%%%%%%%%%%%%%%%%%
Compare these statements with the numerical results depicted in
Figure \ref{fig:2poles1irrat}. Take note that we plot \emph{all} the
zeros of $\Phi_n$ for $n=1,\dots, 150$, in order to reveal the
structure of $Z$ inside $\D$.

Finally, in order to describe the behavior of the zeros that are
closest to the singularity of the weight we must consider again
function $\mathcal H(\beta ; \zeta )$ introduced in
\eqref{def_H_for_local}; let us denote by $h(\beta )$ its zero of
smallest absolute value such that $\Re(h(\beta ))\geq 0$ and
$\Im(h(\beta ))> 0$. According to \eqref{jumpsH_across_imaginary},
$z=-\overline{h(\beta )}$ is also a zero of $\mathcal H(\beta ;
\zeta )$.

\begin{theorem} \label{thm:closestToA}
Among the zeros of $\Phi_n(z)$, the closest to the singularity $a_k$
of the weight are those, given asymptotically by
$$
z^+=a_k\, e^{ 2 i h(\beta_k )/n}  \, \left(1+ o(1)
  \right)\quad \text{and} \quad
z^-=a_k\, e^{ -2 i \, \overline{ h(\beta_k )}/n}  \, \left(1+ o(1)
  \right)\,,
$$
approaching $a_k$ symmetrically with respect to the radius $a_k\cdot
(0,1)$.
\end{theorem}
%%%%%%%%%%%%%%%%%%%%%%%%%%%%%%%%%%%%%%%%%%%%%%%%%%%%%%%%%%
\begin{figure}[htb]
\centering \begin{overpic}[scale=0.83]{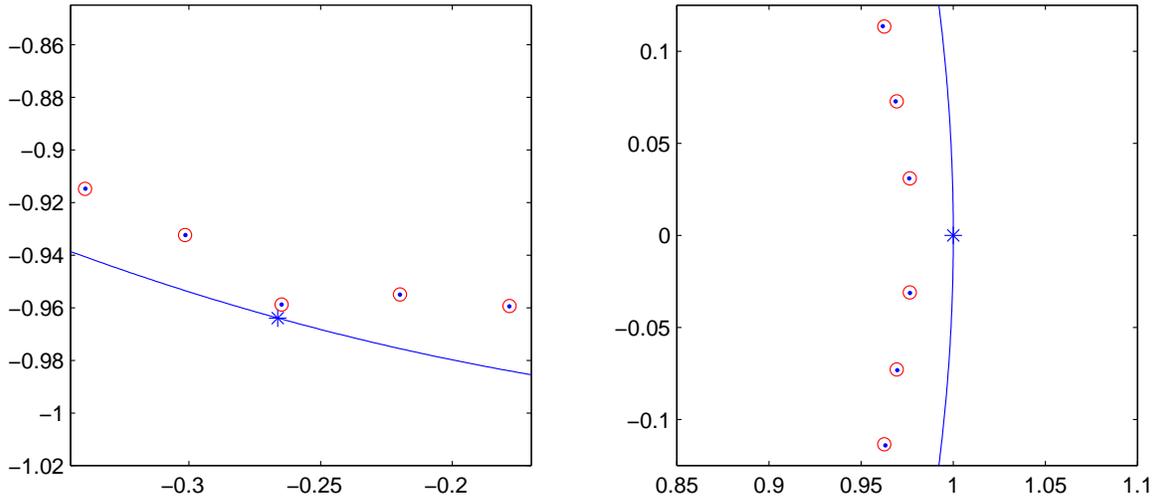}
\end{overpic}
\caption{Zeros of $\Phi_{150}$ (dots) in a neighborhood of $z=a$
(left) and $z=1$ (right) for $W(z)=|z-1|^{1/3} |z-a|^{-2/3}$, $z\in
\T$, where $a=\exp(\pi i \sqrt{2})$. Points $1$ and $a$ are
indicated with asterisks.  Centers of the circles are the zeros of
$\mathcal H(-1/3 ; \zeta_{150}(z) )$ and $\mathcal H(1/6 ;
\zeta_{150}(z) )$, respectively. } \label{fig:local_zeros150}
\end{figure}
%%%%%%%%%%%%%%%%%%%%%%%%%%%%%%%%%%%

\begin{remark} \label{remark:localBehavior:new}
Using formulas (9.6.3) of \cite{abramowitz/stegun:1972} we see that
$\zeta =-i h(\beta)$ must satisfy the following equation:
\begin{equation}\label{combinationOfBessel}
I_{\beta +1/2}(\zeta )-I_{\beta -1/2}(\zeta )=0\,, \quad \Re(\zeta
)>0\,, \quad \Im(\zeta )\geq 0\,,
\end{equation}
where $I_\nu $ is the modified Bessel function. Some facts about the
zeros of the function in the left hand side of
\eqref{combinationOfBessel} have been kindly provided to us by M.\
Muldoon \cite{Muldoon}. For instance, for $\beta>0$ this function is
strictly negative in $(0,+\infty)$ (see \cite{Soni65}) and thus has
no real positive zeros. Furthermore, for $-1/2<\beta <0$ this
function has apparently a unique real positive zero, which is a
monotonically decreasing function of $\beta $. Complementing these
results with numerical experiments (see e.g.\ Remark
\ref{remark:zero} in Section \ref{sec:asymptotics}), we claim that:
\begin{itemize}
\item for $\beta >0$,
$\Re(h(\beta ))> 0$ and $|h(\beta )|>\pi$;

\item for $-1/2<\beta
<0$, $\Re(h(\beta ))= 0$ and $|h(\beta )|<\pi$.
\end{itemize}
In particular, for $-1/2<\beta <0$, $z^+=z^-$, lying (at least,
asymptotically) on the radius $a_k\cdot (0,1)$.
\end{remark}

\begin{remark} \label{remark:localBehavior}
We can compare the rate of convergence of the zeros of $\Phi_n$ to
the unit circle. While, according to Theorem \ref{thm:clock}, the
bulk of the zeros tend to $\T$ with a speed whose leading term is
$n^{-1/n}$ (Theorem \ref{thm:clock}), the zeros closest to the
singularity $a_k$ are attracted by  $\T$ with a rate, roughly
speaking, equal to $c^{-1/n}$, where $c=e^{2 \Im h(\beta _k)}>1$.
Taking into account also that $|z^\pm -a_k|=\frac{2}{n} |h(\beta
)|\, (1+o(1))$, Remark \ref{remark:localBehavior:new} and
\eqref{asymptArgEquidistr}, it explains the influence that $a_k$'s
exert on the neighboring zeros of $\Phi_n$: both pushing them to
$\T$ and ``repelling'' them from (if $\beta
>0$) or attracting them to (if $\beta <0$) the singularity (cf.\ Figure
\ref{fig:2poles1irrat} and zoom in Figure \ref{fig:local_zeros150}).
\end{remark}

\begin{remark}
The method of proof of the statements above can actually handle a
more general situation, when the weight of orthogonality depends on
$n$,
\begin{equation*}%\label{THE_weight_varying}
W_n(z) \isdef  w(z)\, \prod_{k=1}^m |z-a_k|^{2\beta _{k,n}} \,,
\quad z\in \T \,, \quad ( a_k \in \T \text{ for }   k=1,\dots, m)\,,
\end{equation*}
assuming that for values $\beta _k(n)>-1/2$, $k=1, \dots, m$, the
limits
$$
\lim_{n\to \infty} \frac{\beta _{k,n}}{n}=\beta _k\geq 0\,, \quad
k=1,\dots, m\,,
$$
exist.
\end{remark}

The structure of the rest of the paper is as follows. In the next
section we study in more detail the properties of the Szeg\H{o} and
the scattering functions of $W$. The Riemann-Hilbert
characterization and the nonlinear steepest descent methods of Deift
and Zhou (see e.g. \cite{MR2001f:42037,MR2001g:42050}, as well as
the monograph \cite{MR2000g:47048}), are carried out in Section
\ref{sec:RH}, which allows to prove the announced results in Section
\ref{sec:asymptotics}.

\section{Szeg\H{o} and scattering function for $W$}
\label{sec:scattering}

Following the notation introduced above, let $D_{\rm i}(w;z)$ and
$D_{\rm e}(w;z) $ be, respectively, the interior and exterior values
of the Szeg\H{o} function defined by \eqref{standardSzego} for the
analytic component $w$ in \eqref{THE_weight}. Since $w$  is
holomorphic and non-vanishing in $\A_\rho$, then both $D_{\rm
i}(w;z)$ and $D_{\rm e}(w;z) $ admit a holomorphic extension across
$\T $, and maintaining the same notation for these analytic
continuations we have
\begin{equation*}%\label{boundary_value_Szego}
\frac{D_{\rm i}(w; z)}{D_{\rm e}(w; z)}=w(z)\,.
\end{equation*}
This formula gives in a certain way a canonical analytic extension
of $w$ from $\T $ to the annulus $\A_\rho$.

It is convenient to construct explicitly the Szeg\H{o} function for
the modified weight $W$; with this purpose we select the
singe-valued branches of the corresponding functions as follows: the
generalized polynomial
\begin{equation}\label{def_q}
q(z)\isdef \prod_{k=1}^m (z-a_k)^{\beta _k/2}
\end{equation}
is a single-valued analytic function in $\C\setminus
\left(\cup_{k=1}^m a_k\cdot [1,+\infty) \right)$, for which we fix
the value of $q(0)$. Consequently, $\overline{q(1/\bar z)}$ is
single-valued and analytic in  $\C\setminus \left(\cup_{k=1}^m
a_k\cdot [0,1] \right) $, and
$$
\overline{q(1/\bar z)}\, \big|_{z=\infty}= \overline{q(0)}
=\frac{1}{q(0)}.
$$
If we consider each ray $a_k\cdot[1, +\infty)$ oriented towards
infinity, then
\begin{equation}\label{boundary_values_q1}
    q_+(z)=e^{-\pi i \beta _k} q_-(z)\,, \quad z \in a_k\cdot[1,
    +\infty)\,, \quad k=1, \dots, m\,.
\end{equation}
In the same fashion, if $a_k\cdot[0, 1]$ also have the natural
orientation from the origin to $a_k$, then
\begin{equation}\label{boundary_values_q2}
    \left[\overline{q(1/\bar z)}\right]_+=e^{\pi i \beta _k}
    \left[\overline{q(1/\bar z)}\right]_-\,, \quad z \in a_k\cdot(0,
    1)\,, \quad k=1, \dots, m\,.
\end{equation}

With this convention we can write the Szeg\H{o} functions for the
modified weight $W$:
\begin{equation}\label{def:D_eandD_i}
D_{\rm i}(W; z)=   \frac{q^2(z)}{q^2(0)} \, D_{\rm i}(w; z)\,, \quad
D_{\rm e}(W; z)=   \frac{D_{\rm e}(w; z)}{q^2(0)\, \left(
\overline{q (1/\bar z)}\right)^2 }\,.
\end{equation}
Then $D_{\rm i}(W; z)$ is holomorphic in $\D_{1/\rho} \setminus
\left(\cup_{k=1}^m a_k\cdot [1,1/\rho ) \right)$, $D_{\rm e}(W; z)$
is holomorphic in $\{ z\in\C:\, |z|>\rho \} \setminus
\left(\cup_{k=1}^m a_k\cdot (\rho,1] \right)$, and
\begin{equation}\label{tau}
\tau \isdef \frac{1}{D_{\rm i}(W; 0)}=\frac{1}{ D_{\rm i}(w;
0)}=D_{\rm e}(W; \infty)=D_{\rm e}(w; \infty)>0\,.
\end{equation}
Using the definition in \eqref{standardSzego} we see that $\tau $ is
related with the geometric mean defined in \eqref{def_G} by
\begin{equation}\label{tauAlternativeBis}
\tau ^{-2}=G[W]=G[w]\,.
\end{equation}

Furthermore, formula
\begin{equation}\label{boundaryValueforD}
W(z)=q^2(z) \overline{q (1/\bar z)}^2  w(z)=\frac{D_{\rm i}(W;
z)}{D_{\rm e}(W; z)}
\end{equation}
is valid and provides an analytic continuation of the weight $W$ to
the cut annulus $\A_\rho  \setminus \left(\cup_{k=1}^m a_k\cdot
(\rho ,1/\rho ) \right)$. By
\eqref{boundary_values_q1}--\eqref{boundary_values_q2}, with the
orientation of the cuts toward infinity we have
\begin{equation}\label{boundaryValueD_onCuts}
\begin{split}
\left[D_{\rm i}(W; z) \right]_+ & =e^{-2\pi i \beta _k}\left[D_{\rm
i}(W; z) \right]_-\,, \quad z \in a_k\cdot(1,
    1/\rho )\,, \quad k=1, \dots, m\,,\\ \left[D_{\rm e}(W; z) \right]_+ & =e^{-2\pi i \beta _k}\left[D_{\rm
e}(W; z) \right]_-\,, \quad z \in a_k\cdot (\rho ,
    1 )\,, \quad k=1, \dots, m\,.
\end{split}
\end{equation}

Recall that the scattering function for $w$,
\begin{equation*}%\label{def_F_analytic}
    \FF(w;z) = D_{\rm i}(w; z)D_{\rm e}(w; z)\,,
\end{equation*}
is holomorphic in the annulus $\A_\rho $. By \eqref{symmetry},
\begin{equation*}%\label{symmetryOfF}
 \overline{\FF\left(w;\frac{1}{\overline{z}}
\right)}=\frac{1}{\FF(w;z)}\,, \quad \text{for } z\in \A_\rho \,,
\qquad
 |\FF(w;z)|=1 \text{ on } \T \,.
\end{equation*}
With the definition \eqref{def_Scattering} and formulas
\eqref{def:D_eandD_i} we have
\begin{equation*}%\label{def_F_branch}
    \FF(W;z) = D_{\rm i}(W; z)D_{\rm e}(W; z)=
  \left( \frac{  q (z) }{q^2(0)\, \overline{q (1/\bar
z)}  } \right)^2 \FF(w; z)   \,,
\end{equation*}
that is also analytic and single-valued in the cut annulus $\A_\rho
\setminus \left(\cup_{k=1}^m a_k\cdot (\rho ,1/\rho ) \right)$. By
\eqref{boundaryValueD_onCuts},
\begin{equation}\label{boundaryValueF_onCuts}
\FF_+(W; z)  =e^{-2\pi i \beta _k} \FF_-(W; z) \,, \quad z \in
a_k\cdot(\rho ,
    1/\rho )\setminus \{a_k \}\,, \quad k=1, \dots, m\,.
\end{equation}
It is straightforward also to check that $\FF(W;z)$ is bounded in a
neighborhood of each $a_k$.

Recall the definition in \eqref{notationBC}. With our assumptions on
$\delta>0 $ the closed disks $\BB_k\cup \CC_k$  are disjoint, and
$\FF(w;z)$ is analytic in each $\BB_k$. Furthermore, motivated by
\eqref{boundaryValueF_onCuts} we define
\begin{equation}\label{defSModified}
    \widehat \FF_k(W;z)\isdef \begin{cases}
e^{ \pi i \beta _k} \FF(W; z)\,, & \text{if } z\in \BB_k \text{ and
} \arg(z)>\arg(a_k)\,, \\ e^{- \pi i \beta _k} \FF(W; z)\,, &
\text{if } z\in \BB_k \text{ and } \arg(z)< \arg(a_k)\,,
    \end{cases} \quad k=1, \dots, m\,.
\end{equation}
Then $\widehat \FF_k(W;z)$ is holomorphic in $\BB_k$, $k=1, \dots,
m$, and we define
\begin{equation}\label{def_new_constant}
    \vartheta_k \isdef \widehat \FF_k(W;a_k)\in \T \,, \quad k=1, \dots,
    m\,.
\end{equation}
These constants have been used in the formulation of several results
in the previous section.
\begin{remark}
Formally
$$
 \vartheta_k=\prod_{j\neq k} \left( -\frac{a_k}{a_j}\right)^{\beta
 _k} \FF(w;a_k)\,,
$$
but the selection of the right branch in each case should be
specified. So the definition by \eqref{def_new_constant} is
preferred in order to avoid ambiguity.
\end{remark}

\section{Riemann-Hilbert analysis for orthogonal polynomials}
\label{sec:RH}

We assume the unit circle $\T $ oriented counterclockwise. The
starting point of all the analysis is the fact that under
assumptions  above conditions \eqref{orthogonalityConditions} can be
rewritten in terms of a non-hermitian orthogonality for $\varphi_n$
and $\varphi_n^*$:
\begin{align*}
 \oint_{\T } \varphi_n(z)  z^{n-k-1} \,
\frac{W(z)}{z^n}\, dz & =0,\quad \text{for } k=0, 1, \dots, n-1\,,
 \\
 \oint_{\T } \varphi^*_{n-1}(z)  z^{k} \,
\frac{W(z)}{z^n} \, dz & = \begin{cases} 0, & k=0, 1, \dots, n-2, \\
i/\kappa _{n-1}, & k=n-1\,.
\end{cases}
\end{align*}
By standard arguments (see e.g.\ \cite{MR2000e:05006} or
\cite{MR2000g:47048}, as well as the seminal paper \cite{Fokas92}
where the Riemann-Hilbert approach to orthogonal polynomials
started),
\begin{equation}\label{Y}
Y(z)=\begin{pmatrix} \Phi_n(z) & \displaystyle \dfrac{1}{2\pi
i}\,\oint_{\strut\T }
\dfrac{\Phi_n(t) W(t)\, dt}{t^n(t-z)}  \\
-2\pi \kappa _{n-1} \varphi^*_{n-1}(z) & - \, \displaystyle
\dfrac{\kappa _{n-1}}{ i}\,\oint_{\strut\T }
\dfrac{\varphi^*_{n-1}(t) W(t)\, dt}{t^n(t-z)}
\end{pmatrix}
\end{equation}
is a unique solution of the following Riemann-Hilbert problem: $Y$
is holomorphic in $\C\setminus {\T }$,
\begin{equation*}%\label{RHproblem}
Y_+(z)=Y_-(z)\, \begin{pmatrix} 1 & W(z)/z^n \\ 0 & 1
\end{pmatrix}\,, \quad z\in \T \,, \qquad %\text{and} \quad
\lim_{z\to \infty} Y(z)\,\begin{pmatrix} z^{-n} & 0 \\ 0 & z^n
\end{pmatrix}=I\,,
\end{equation*}
where $I$ is the $2 \times 2$ identity matrix, and
$$
Y(z)=\begin{cases} \mathcal O \begin{pmatrix}  1 & |z-a_k|^{2\beta
_k} \\ 1 & |z-a_k|^{2\beta _k}
\end{pmatrix} \,, & \text{if } \beta _k<0 \,, \\ \mathcal O \begin{pmatrix}  1 & 1 \\ 1 &
1 \end{pmatrix} \,, & \text{if } \beta _k\geq 0\,,
\end{cases}
$$
as $z\to a_k$, $z\in \C\setminus \T $, and $k=1, \dots, m$.

In order to perform the steepest descent analysis as described in
\cite{MR2000g:47048} (see also \cite{Kuijlaars:03}) we build a
series of \emph{explicit} and \emph{reversible} steps in order to
arrive at an equivalent problem, which is solvable, at least in an
asymptotic
sense. We will use the following notation: $\sigma_3 =\begin{pmatrix} 1 & 0 \\
0& -1 \end{pmatrix}$ is the Pauli matrix, and for any non-zero
$x$ and integer $m$, $ x^{\sigma_3}=\begin{pmatrix} x & 0 \\
0& 1/x \end{pmatrix} $.

\subsection{Global analysis}

If we define
\begin{equation}\label{defH}
H(z)\isdef \begin{cases} z^{-n \sigma _3} , & \text{if } |z|>1, \\
I, & \text{if } |z|<1,
\end{cases}
\end{equation}
and put $ T(z) \isdef Y(z)\, H(z) $, then $T$ becomes holomorphic in
$\C\setminus {\T  }$, with $ \lim_{z\to \infty} T(z)=I $, and
satisfying the jump condition
\begin{equation*}%\label{jumpT}
T_+(z)=T_-(z)\, \begin{pmatrix} z^n & W(z)  \\ 0 & z^{-n}
\end{pmatrix}
\end{equation*}
on ${\T}$. Obviously, $T$ exhibits the same local behavior at $\AA$
as $Y$.

\begin{figure}[htb]
\centering \begin{overpic}[scale=0.7]{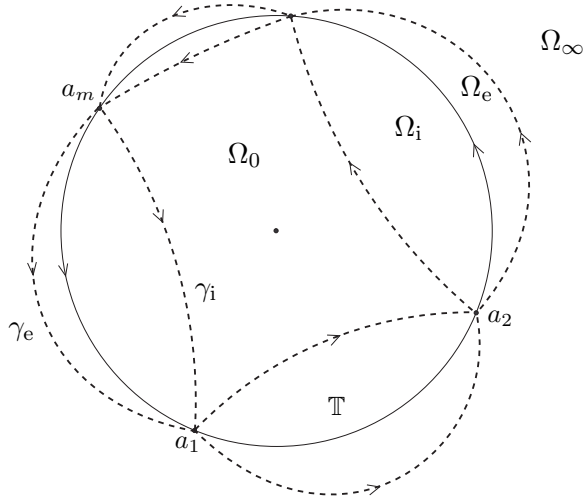}
      \put(55,27){\small $\T $}
         \put(35,45){$\gamma_{\rm i}$}
         \put(7,39){$\gamma_{\rm e}$}
          \put(40,65){$\Omega_0$}
          \put(65,69){$\Omega_{\rm i}$}
          \put(75,75){$\Omega_{\rm e}$}
          \put(87,82){$\Omega_\infty$}
          \put(15,75){\small $a_m$}
          \put(32,22){\small $a_1$}
         \put(79,41){\small $a_2$}
\end{overpic}
\caption{Opening lenses.}\label{fig:zeros_case1}
\end{figure}

Let $\gamma_{\rm i}$ be a closed Jordan contour, piecewise analytic,
entirely contained in the cut annulus $\{z\in \C:\, \rho   <|z|<1\}
\setminus \cup_{k=1}^m a_k \cdot (\rho ,1)$, except for points
$a_k$: $ \gamma_{\rm i}\cap \T =\AA$. Let also $\gamma_{\rm e}\isdef
\{ 1/\bar z: z\in \gamma_{\rm i}\}$ the contour symmetric to $\gamma
_{\rm i}$ with respect to $\T $. We take both $\gamma _{\rm i}$ and
$\gamma _{\rm e}$ oriented counterclockwise. Let $\Omega_0$ be the
connected component of $\C\setminus \gamma_{\rm i}$ containing the
origin, and $\Omega_\infty$ the corresponding unbounded component of
$\C\setminus \gamma_{\rm e}$. Furthermore, we denote by
$$
\Omega_{\rm i}\isdef \D \setminus \overline{\Omega_0}\,, \qquad
\Omega_{\rm e}\isdef \{ 1/\bar z: z\in \Omega_{\rm i}\}
$$
(see Fig. \ref{fig:zeros_case1}).

%%%%%%%%%%%%%%%%%%%%%%%%%%%%%%%%%%%%%%%%%%%%%%%%%%%%%%%%%%
%%%%%%%%%%%%%%%%%%%%%%%%%%%%%%%%%%%%

Then $W$ is correctly defined in both $\Omega_{\rm i}$ and
$\Omega_{\rm e}$, and we can set
\begin{equation}\label{KforT1}
    K(z)\isdef \begin{cases} I, & \text{if } z \in \Omega_0 \cup \Omega_\infty, \\
\begin{pmatrix} 1 & 0  \\ z^{n}/  W(z)  & 1
\end{pmatrix}^{-1}\,, & \text{if }
z\in \Omega_{\rm i}, \\ \begin{pmatrix}  1 & 0  \\ 1/(z^{n } W(z)) &
1
\end{pmatrix}\,, & \text{if }
z\in \Omega_{\rm e}.
\end{cases}
\end{equation}
Using $K$ we make a new transformation: $U(z)\isdef T(z)K(z)$.
 Matrix valued function $U$ is holomorphic in $\C\setminus \left( \T  \cup \gamma_{\rm i} \cup
\gamma_{\rm e}\right)$,  $ \lim_{z\to \infty} U(z)=I$, and
$$
U_+(z)=U_-(z)\, J_U(z), \quad z\in  \T  \cup \gamma_{\rm i} \cup
\gamma_{\rm e},
$$
where
\begin{equation}\label{defJumpU1}
J_U(z)=\begin{cases} \begin{pmatrix} 0 &  W(z)  \\ -1/  W(z)  & 0
\end{pmatrix}, & \text{if } z\in \T  \setminus \AA, \\   \begin{pmatrix} 1 & 0  \\ z^{n }/  W(z)  & 1
\end{pmatrix}, & \text{if } z\in \gamma_{\rm i}\setminus \AA, \\
\begin{pmatrix} 1 & 0  \\ 1/(z^{n } W(z)) & 1
\end{pmatrix}, & \text{if } z\in \gamma_{\rm e}\setminus \AA.
\end{cases}
\end{equation}
Moreover, the local behavior for $U(z)$ as $z \to a_k$ from $
\Omega_{\rm e} \cup \Omega_{\rm i}$ is now
\begin{equation*}%\label{localBehaviorU}
U(z)=\begin{cases} \mathcal O \begin{pmatrix} 1 & |z-a_k|^{ 2\beta_k
} \\ 1 & |z-a_k|^{ 2\beta_k }
\end{pmatrix}, & \text{if } \beta _k< 0\,, \\
\mathcal O \begin{pmatrix} |z-a_k|^{-2\beta_k } & 1 \\
|z-a_k|^{-2\beta_k } & 1
\end{pmatrix}, & \text{if } \beta _k\geq 0\,.
\end{cases}
\end{equation*}

 With $D_{\rm e}(W; \cdot)$ and
$D_{\rm i}(W;  \cdot)$ defined in \eqref{def:D_eandD_i} and constant
$\tau $ introduced in \eqref{tau}, let
\begin{equation}\label{equ:defN}
    N(z)\isdef \begin{cases}
   \left( \dfrac{D_{\rm e}(W; z)}{\tau} \right)^{\sigma _3} , & \text{if } |z|>1\,, \\
  \begin{pmatrix} 0 & D_{\rm i}(W; z)/\tau   \\ -\tau /D_{\rm i}(W; z) & 0
\end{pmatrix}= \left( \dfrac{D_{\rm i}(W; z)}{\tau} \right)^{\sigma _3} \,
\begin{pmatrix} 0 & 1   \\ -1 & 0
\end{pmatrix} , & \text{if } |z|<1\,.
    \end{cases}
\end{equation}
By \eqref{boundaryValueforD}, $N$ has the same jump on $\T $ as $U$,
and by \eqref{tau}, it exhibits the same behavior at infinity.
Hence, $U(z)N^{-1}(z)$ tends to $I$ as $z\to \infty$, is holomorphic
in $\C\setminus (\T  \cup \gamma_{\rm e}\cup \gamma_{\rm i})$, and
has jumps across these curves asymptotically close to $I$, except
when we approach the singular set $\AA$. We have to handle the
behavior at each individual singular point $a_k$ by means of the
local analysis.

\subsection{Local analysis}

Let us pick a singular point $a_\ell \in \AA$. For the sake of
brevity along this subsection we use the following shortcuts for the
notation: $a\isdef a_\ell $, $\beta \isdef \beta _\ell $, $\BB\isdef
\BB_\ell $, $\CC \isdef \CC_\ell $ (where $\BB_\ell $, $\CC_\ell $
and $\delta $ were defined in \eqref{notationBC}), and $\BB^+\isdef
\{z\in \BB:\, \arg(z)>\arg(a) \}$, $\BB^-\isdef \{z\in \BB:\,
\arg(z)<\arg(a) \}$. We also write $\widehat \Omega_j \isdef
\Omega_j \cap \BB$, where $j \in \{{\rm i}, {\rm e} , 0, \infty \}$,
and analogous notation for curves: $\widehat \T \isdef \T  \cap
\BB$, etc. Finally, $ \widehat \FF (W;z)\isdef  \widehat
\FF_\ell(W;z) $ (see the definition in \eqref{defSModified}).

Our goal is to build a matrix $P(z)\isdef P(a,\beta ;z)$ meeting the
following requirements:
\begin{enumerate}
\item[(P1)] $P$ is holomorphic in $\BB \setminus (\T  \cup \gamma_{\rm
i} \cup \gamma_{\rm e})$ and satisfies across $\widehat \T  \cup
\widehat \gamma_{\rm i} \cup \widehat \gamma_{\rm e}$ the jump
relation $P_+(t)=P_-(t) J_U(t)$, with $J_U$ defined in
\eqref{defJumpU1}.

\item[(P2)] $P(z)$ has the following local behavior as $z \to a$: if $\beta
<0$, then
\begin{equation*}%\label{localBehaviorPbetaless0}
P(z)=  \mathcal O \begin{pmatrix}  1 & |z-a|^{2\beta }
\\ 1 & |z-a|^{2\beta }
\end{pmatrix} \,,
\end{equation*}
and if $\beta \geq 0$,
\begin{equation*}%\label{localBehaviorPbetamore0}
P(z)=\begin{cases} \mathcal O \begin{pmatrix} 1 & 1 \\ 1 & 1
\end{pmatrix}, & \text{from } \widehat \Omega_{0} \cup \widehat
\Omega_{\infty} \,, \\
\mathcal O \begin{pmatrix} |z-a|^{-2\beta  } & 1 \\
|z-a|^{-2\beta  } & 1
\end{pmatrix}, & \text{from } \widehat \Omega_{\rm e} \cup \widehat \Omega_{\rm i}\,.
\end{cases}
\end{equation*}

\item[(P3)] $P(z)$ matches $N(z)$ on $\CC$, in the sense $P(z)
N^{-1}(z)=I+\mathcal O (n^{-1})$ for $z\in \CC$.
\end{enumerate}

%%%%%%%%%%%%%%%%%%%%%%%%%%%%%%%%%%%%%%%%%%%%%%%%%%%%%%%%%%
\begin{figure}[htb]
\centering \begin{overpic}[scale=0.9]{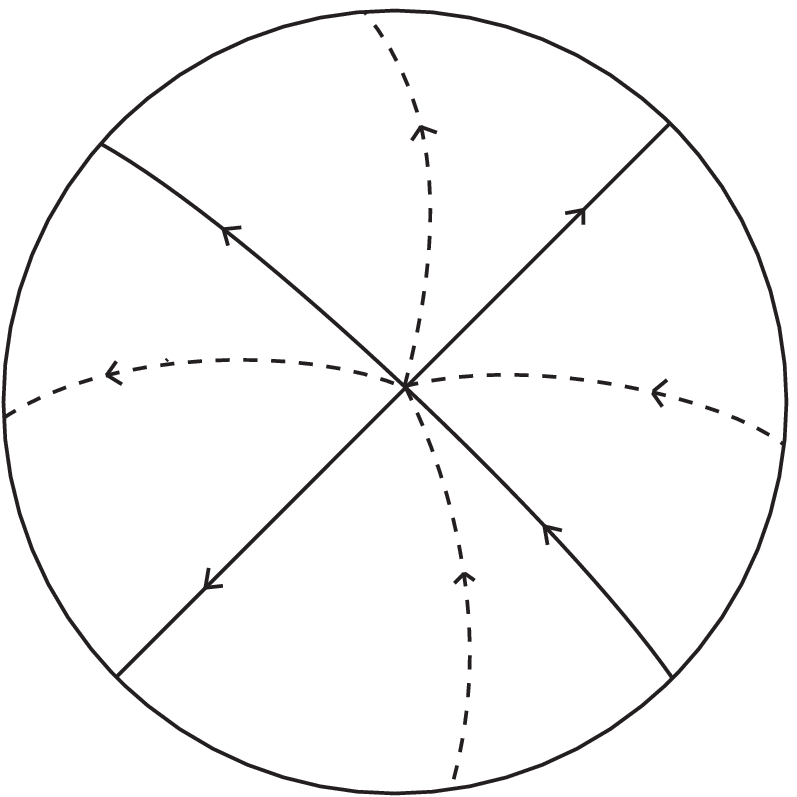}
         \put(50,30){$\gamma_{\rm i}$}
         \put(20,50){$\gamma_{\rm i}$}
         \put(48,75){$\gamma_{\rm e}$}
         \put(85,45){$\gamma_{\rm e}$}
          \put(40,20){$\Omega_0$}
          \put(20,40){$\Omega_0$}
           \put(65,85){$\Omega_\infty$}
          \put(80,65){$\Omega_\infty$}
          \put(35,85){$\Omega_{\rm e}$}
          \put(85,35){$\Omega_{\rm e}$}
          \put(12,62){$\Omega_{\rm i}$}
          \put(70,15){$\Omega_{\rm i}$}
                    \put(17,25){$\Gamma_{\rm i}$}
                   \put(80,75){$\Gamma_{\rm e}$}
         \put(49,47){\small $a$}
       \put(15,73){\small $\T $}
\end{overpic}
\caption{Local analysis in $\BB$.}\label{fig:local1}
\end{figure}
%%%%%%%%%%%%%%%%%%%%%%%%%%%%%%%%%%%%

As a first step  we  reduce the problem to the one with constant
jumps. Let us denote $\Gamma_{\rm i}\isdef a\cdot (1-\delta, 1)$ and
$\Gamma_{\rm e}\isdef a\cdot (1, 1+\delta)$, oriented both from
$z=a$ to infinity (see Fig. \ref{fig:local1}). Let $w^{1/2}(z) $ and
$z^{1/2}$ denote the principal holomorphic branches of these
functions in $b$, and (cf.\ \eqref{boundaryValueforD}),
\begin{equation*}%\label{def_sqrtW}
W^{1/2}(z)\isdef q(z) \overline{q(1/\overline{z})}\, w^{1/2}(z)
\end{equation*}
with $q$ defined in \eqref{def_q}. Then $W^{1/2}$ is holomorphic in
$\BB \setminus a\cdot (1-\delta,1+\delta)$, and according to
\eqref{boundary_values_q1}--\eqref{boundary_values_q2},
\begin{equation*}%\label{jumpOfw}
\frac{W^{1/2}_+(z)}{W^{1/2}_-(z)}= e^{- \pi i \beta }\quad \text{on
} \Gamma_{\rm i}\,, \quad \text{and} \quad
\frac{W^{1/2}_+(z)}{W^{1/2}_-(z)}= e^{  \pi i \beta }\quad \text{on
} \Gamma_{\rm e}\,.
\end{equation*}
Thus, if we define
\begin{equation}\label{lambda}
\lambda( \beta ;z) \isdef \begin{cases} e^{\pi i \beta } W^{1/2}(z)
z^{n/2}, & z
\in  (\widehat \Omega_{\rm e}  \cup \widehat \Omega_\infty   ) \cap \BB^+  ,\\
e^{-\pi i \beta } W^{1/2}(z) z^{n/2}, & z \in  (\widehat \Omega_{\rm
e}  \cup \widehat \Omega_\infty   ) \cap \BB^-  ,\\  e^{-\pi i \beta
} W^{1/2}(z) z^{-n/2}, & z
\in  (\widehat  \Omega_{\rm i}  \cup \widehat \Omega_0   ) \cap \BB^+  ,\\
e^{\pi i \beta } W^{1/2}(z) z^{-n/2}, & z \in  (\widehat \Omega_{\rm
i}  \cup \widehat \Omega_0   ) \cap \BB^-  ,
\end{cases}\,,
\end{equation}
and set
\begin{equation}\label{def_R}
R(z) \isdef  P(z)\, \lambda (  \beta ; z)^{  \sigma _3} \,, \quad z
\in b \setminus (\Gamma_{\rm i} \cup \Gamma_{\rm e} \cup \T  \cup
\gamma_{\rm i} \cup \gamma_{\rm e})\,,
\end{equation}
we get for $R$ the following problem: $R$ is holomorphic in $b
\setminus (\Gamma_{\rm i} \cup \Gamma_{\rm e} \cup \T  \cup
\gamma_{\rm i} \cup \gamma_{\rm e})$, and satisfies the jump
relation $R_+(z)=R_-(z) J_R(z)$, with
\begin{equation*}%\label{defJumpR1}
J_R(z)=\begin{cases}  \begin{pmatrix} 0 & 1  \\ -1  & 0
\end{pmatrix}, & \text{if } z \in \widehat \T   , \\
 \begin{pmatrix} 1 & 0  \\  e^{-2\pi i \beta  }  & 1
\end{pmatrix}, & \text{if } z\in   (\widehat \gamma_{\rm i}   \cap \BB^+)\cup
(\widehat \gamma_{\rm e}   \cap \BB^-)  ,
\\
\begin{pmatrix} 1 & 0  \\  e^{ 2\pi i \beta  }  & 1
\end{pmatrix}, & \text{if } z\in   (\widehat \gamma_{\rm i}   \cap \BB^-)
\cup (\widehat \gamma_{\rm e}   \cap \BB^+) ,\\
\begin{pmatrix} e^{ \pi i \beta } & 0  \\  0  & e^{-\pi i \beta }
\end{pmatrix}, & \text{if } z\in \Gamma_{\rm i}\cup \Gamma_{\rm e}\,.
\end{cases}
\end{equation*}
Moreover, $R$ has the following local behavior as $z\to a$: if
$\beta \geq 0$,
\begin{equation*}%\label{behavior_R1}
R(z)=\begin{cases} \mathcal O \begin{pmatrix}  |z-a|^{\beta  } & |z-a|^{-\beta  } \\
|z-a|^{\beta  } &
|z-a|^{-\beta  } \end{pmatrix}, & \text{if } z \in \widehat \Omega_0 \cup \widehat \Omega_\infty, \\
\mathcal O \begin{pmatrix}  |z-a|^{-\beta  } & |z-a|^{-\beta  } \\
|z-a|^{-\beta  } & |z-a|^{-\beta  }
\end{pmatrix}, & \text{if } z \in \widehat \Omega_{\rm e} \cup \widehat \Omega_{\rm i},
\end{cases}
\end{equation*}
and if $\beta <0$, then
\begin{equation*}%\label{behavior_R2}
R(z)=   \mathcal O \begin{pmatrix}  |z-a|^{\beta  } & |z-a|^{\beta  } \\
|z-a|^{\beta  } & |z-a|^{\beta  } \end{pmatrix}.
\end{equation*}

As it could be expected, this problem resembles very much the one we
face during the local analysis for the generalized Jacobi weight on
the real line (see \cite[Theorem 4.2]{Vanlessen03}, \cite{MR2087231}
and \cite{Kuijlaars/Vanlessen:03}). We take advantage of the results
proved therein in order to abbreviate the exposition.

Let us define $\Sigma_i\isdef \{x e^{i\pi /4}:\, x\geq 0 \}$, $i=1,
\dots, 8$; we take $\Sigma_i$ oriented towards the origin for $i=3,
4, 5$, and the rest towards infinity (see Fig. \ref{fig:local2}).
This splits the plane into eight regions, as indicated there, marked
with Roman numbers $I$ to $VIII$. Using Hankel functions
$H_\nu^{(i)}$ and modified Bessel functions $K_\nu $ and $I_\nu$ we
build a piece-wise analytic matrix-valued function
$\Psi=\Psi(\beta;\cdot) $, $\beta
>-1/2$, in the following manner (where $\zeta^{1/2}$ denotes the main branch in $(-\infty,
0]$):

%%%%%%%%%%%%%%%%%%%%%%%%%%%%%%%%%%%%%%%%%%%%%%%%%%%%%%%%%%
\begin{figure}[htb]
\centering
\begin{overpic}[scale=0.8]{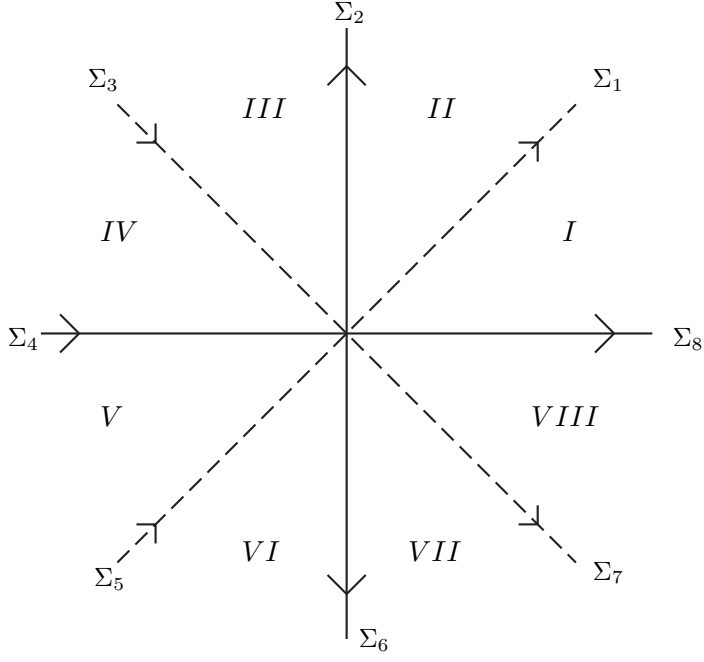}
\put(85,65){$I$} \put(63,85){$II$} \put(33,85){$III$}
\put(10,65){$IV$} \put(10,35){$V$} \put(33,13){$VI$}
\put(60,13){$VII$} \put(80,35){$VIII$}
 \put(90,90){\small $\Sigma_1$}
 \put(48,101){\small $\Sigma_2$}
 \put(8,90){\small $\Sigma_3$}
 \put(-5,48){\small $\Sigma_4$}
 \put(9,9){\small $\Sigma_5$}
 \put(52,-1){\small $\Sigma_6$}
 \put(90,10){\small $\Sigma_7$}
 \put(103,48){\small $\Sigma_8$}
\end{overpic}
\caption{Local parametrix.}\label{fig:local2}
\end{figure}
%%%%%%%%%%%%%%%%%%%%%%%%%%%%%%%%%%%%%%%%%%%%%%%%%%%%%%%%%%

For $\zeta\in$ I,
\begin{equation}\label{RHPPSIx0solution1}
    \Psi  (\beta ; \zeta ) \isdef \frac{1}{2}\sqrt\pi\zeta^{1/2}
    \begin{pmatrix}
        H_{\beta  +\frac{1}{2}}^{(2)}(\zeta) &
            -i H_{\beta  +\frac{1}{2}}^{(1)}(\zeta) \\[2ex]
        H_{\beta  -\frac{1}{2}}^{(2)}(\zeta) &
            -i H_{\beta  -\frac{1}{2}}^{(1)}(\zeta)
    \end{pmatrix}
    e^{-(\beta  +\frac{1}{4}) \pi i \sigma_3}.
\end{equation}
For $\zeta\in$ II,
\begin{equation}\label{RHPPSIx0solution2}
    \Psi  (\beta ; \zeta ) \isdef
    \begin{pmatrix}
        \sqrt\pi\zeta^{1/2}I_{\beta  +\frac{1}{2}}(\zeta e^{-\frac{\pi i}{2}}) &
            -\frac{1}{\sqrt\pi}\zeta^{1/2}K_{\beta  +\frac{1}{2}}(\zeta e^{-\frac{\pi i}{2}}) \\[2ex]
        -i\sqrt\pi\zeta^{1/2}I_{\beta  -\frac{1}{2}}(\zeta e^{-\frac{\pi i}{2}}) &
            -\frac{i}{\sqrt\pi}\zeta^{1/2}K_{\beta  -\frac{1}{2}}(\zeta e^{-\frac{\pi i}{2}})
    \end{pmatrix}
    e^{-\frac{1}{2}\beta  \pi i\sigma_3}.
\end{equation}
For $\zeta\in$ III,
\begin{equation}\label{RHPPSIx0solution3}
    \Psi  (\beta ; \zeta ) \isdef
    \begin{pmatrix}
        \sqrt\pi\zeta^{1/2}I_{\beta  +\frac{1}{2}}(\zeta e^{-\frac{\pi i}{2}}) &
            -\frac{1}{\sqrt\pi}\zeta^{1/2}K_{\beta  +\frac{1}{2}}(\zeta e^{-\frac{\pi i}{2}}) \\[2ex]
        -i\sqrt\pi\zeta^{1/2}I_{\beta  -\frac{1}{2}}(\zeta e^{-\frac{\pi i}{2}}) &
            -\frac{i}{\sqrt\pi}\zeta^{1/2}K_{\beta  -\frac{1}{2}}(\zeta e^{-\frac{\pi i}{2}})
    \end{pmatrix}
    e^{\frac{1}{2}\beta  \pi i\sigma_3}.
\end{equation}
For $\zeta\in$ IV,
\begin{equation}\label{RHPPSIx0solution4}
    \Psi  (\beta ; \zeta ) \isdef \frac{1}{2}\sqrt\pi(-\zeta)^{1/2}
    \begin{pmatrix}
        i H_{\beta  +\frac{1}{2}}^{(1)}(-\zeta) &
            - H_{\beta  +\frac{1}{2}}^{(2)}(-\zeta) \\[2ex]
        -i H_{\beta  -\frac{1}{2}}^{(1)}(-\zeta) &
            H_{\beta  -\frac{1}{2}}^{(2)}(-\zeta)
    \end{pmatrix}
    e^{(\beta  +\frac{1}{4})\pi i \sigma_3}.
\end{equation}
For $\zeta\in$ V,
\begin{equation}\label{RHPPSIx0solution5}
    \Psi  (\beta ; \zeta ) \isdef \frac{1}{2}\sqrt\pi(-\zeta)^{1/2}
    \begin{pmatrix}
        - H_{\beta  +\frac{1}{2}}^{(2)}(-\zeta) &
            -i H_{\beta  +\frac{1}{2}}^{(1)}(-\zeta) \\[2ex]
        H_{\beta  -\frac{1}{2}}^{(2)}(-\zeta) &
            i H_{\beta  -\frac{1}{2}}^{(1)}(-\zeta)
    \end{pmatrix}
    e^{-(\beta  +\frac{1}{4})\pi i \sigma_3}.
\end{equation}
For $\zeta\in$ VI,
\begin{equation}\label{RHPPSIx0solution6}
    \Psi  (\beta ; \zeta ) \isdef
    \begin{pmatrix}
        -i\sqrt{\pi}\zeta^{1/2}I_{\beta  +\frac{1}{2}}(\zeta e^{\frac{\pi i}{2}}) &
            -\frac{i}{\sqrt\pi}\zeta^{1/2}K_{\beta  +\frac{1}{2}}(\zeta e^{\frac{\pi i}{2}}) \\[2ex]
        \sqrt{\pi}\zeta^{1/2}I_{\beta  -\frac{1}{2}}(\zeta e^{\frac{\pi i}{2}}) &
            -\frac{1}{\sqrt\pi}\zeta^{1/2}K_{\beta  -\frac{1}{2}}(\zeta e^{\frac{\pi i}{2}})
    \end{pmatrix}
    e^{-\frac{1}{2}\beta  \pi i \sigma_3}.
\end{equation}
For $\zeta\in$ VII,
\begin{equation}\label{RHPPSIx0solution7}
    \Psi  (\beta ; \zeta ) \isdef
    \begin{pmatrix}
        -i\sqrt{\pi}\zeta^{1/2}I_{\beta  +\frac{1}{2}}(\zeta e^{\frac{\pi i}{2}}) &
            -\frac{i}{\sqrt\pi}\zeta^{1/2}K_{\beta  +\frac{1}{2}}(\zeta e^{\frac{\pi i}{2}}) \\[2ex]
        \sqrt{\pi}\zeta^{1/2}I_{\beta  -\frac{1}{2}}(\zeta e^{\frac{\pi i}{2}}) &
            -\frac{1}{\sqrt\pi}\zeta^{1/2}K_{\beta  -\frac{1}{2}}(\zeta e^{\frac{\pi i}{2}})
    \end{pmatrix}
    e^{\frac{1}{2}\beta  \pi i \sigma_3}.
\end{equation}
And finally, for $\zeta\in$ VIII,
\begin{equation}\label{RHPPSIx0solution8}
    \Psi  (\beta ; \zeta ) \isdef  \frac{1}{2}\sqrt\pi\zeta^{1/2}
    \begin{pmatrix}
        -i H_{\beta  +\frac{1}{2}}^{(1)}(\zeta) &
            -H_{\beta  +\frac{1}{2}}^{(2)}(\zeta) \\[2ex]
        -iH_{\beta  -\frac{1}{2}}^{(1)}(\zeta) &
            -H_{\beta  -\frac{1}{2}}^{(2)}(\zeta)
    \end{pmatrix}
    e^{(\beta  +\frac{1}{4}) \pi i \sigma_3}.
\end{equation}

\begin{proposition}[\cite{Vanlessen03}, Theorem 4.2] \label{lemma:Vanlessen}
Function $\Psi=\Psi(\beta ; \cdot) $ defined above is holomorphic in
$\C \setminus \bigcup_{i=1}^8 \Sigma_i$, and exhibits the following
jumps:
\begin{align}
 \label{JumpsForPsi1}   \Psi_+(\zeta) & =\Psi_-(\zeta)\, \begin{pmatrix} 0 & 1  \\ -1  & 0
\end{pmatrix}, \quad \text{for } \zeta\in \Sigma_4\cup \Sigma_8,
\\
\Psi_+(\zeta) & =\Psi_-(\zeta)\,   \begin{pmatrix} 1 & 0  \\
e^{-2\pi i \beta  }  & 1
\end{pmatrix}, \quad \text{for } \zeta\in \Sigma_1\cup \Sigma_5, \\
\Psi_+(\zeta) & =\Psi_-(\zeta)\,   \begin{pmatrix} 1 & 0  \\
e^{2\pi i \beta  }  & 1
\end{pmatrix}, \quad \text{for } \zeta\in \Sigma_3\cup \Sigma_7,
\\
\Psi_+(\zeta) & =\Psi_-(\zeta)\,  \begin{pmatrix} e^{ \pi i \beta }
& 0  \\  0  & e^{-\pi i \beta }
\end{pmatrix}, \quad \text{for }  \zeta \in \Sigma_2 \cup \Sigma_6. \label{JumpsForPsi4}
\end{align}
Additionally, if $\beta  \geq 0$ then as $\zeta\to 0$,
\begin{equation*}%\label{asymptBehavior}
\Psi(\zeta)=\begin{cases} \mathcal O \begin{pmatrix}  |\zeta|^{\beta  } & |\zeta|^{-\beta  } \\
|\zeta|^{\beta  } & |\zeta|^{-\beta  } \end{pmatrix}, & \text{if $
\zeta$   is in domains II, III, VI, VII},  \\
\mathcal O \begin{pmatrix}  |\zeta|^{-\beta  } & |\zeta|^{-\beta  } \\
|\zeta|^{-\beta  } & |\zeta|^{-\beta  }
\end{pmatrix}, & \text{otherwise},
\end{cases}
\end{equation*}
and if $\beta  < 0$ then
\begin{equation*}%\label{asymptBehaviorBis}
\Psi(\zeta)=  \mathcal O \begin{pmatrix}  |\zeta|^{\beta  } & |\zeta|^{ \beta  } \\
|\zeta|^{\beta  } & |\zeta|^{ \beta  } \end{pmatrix}.
\end{equation*}
\end{proposition}

Consider in $\C \setminus (-\infty, 0)$ the transformation
\begin{equation*}%\label{def_Transf_F}
 f(z)=-i\, \log(z/a) \,,
\end{equation*}
where we take the principal branch of the logarithm. Then $f$ is a
conformal 1-1 map of $\BB$ onto a neighborhood of the origin.
Moreover, $\T $ is mapped onto $\R$ oriented positively, and we may
use the freedom in the selection of the contours deforming them in
such a way that $f(\widehat \gamma_{\rm i})$ and $f(\widehat
\gamma_{\rm e})$ follow the rays $\Sigma_i$ in $\BB$ with odd
indices $i$ (dashed lines, see Figs.\ \ref{fig:local1} and
\ref{fig:local2}). By construction, matrix
\begin{equation*}%\label{def_Q}
    \Psi \left( \beta ;  \frac{n}{2}\, f(z)  \right) %\, \lambda (  \beta ; z)^{ - \sigma _3}
\end{equation*}
matches the jumps and the local behavior of $R$ in $\BB$. Since a
left multiplication by a holomorphic function has no influence on
the jumps, and taking into account \eqref{def_R}, we see that matrix
$P$ can be built of the form
\begin{equation}\label{formOfP}
P(z) = E(z)\, \Psi \left( \beta ;  \frac{n}{2}\, f(z)  \right)
\lambda (  \beta ; z)^{ - \sigma _3}\,,
\end{equation}
where $E$ is any holomorphic function in $\BB$. An adequate
selection of $E$ is motivated by the matching requirement $P(z)
N^{-1}(z )=I+\mathcal O (n^{-1})$ on the boundary $\CC$.

For the sake of brevity let us denote
\begin{equation}\label{def_zeta_local_analysis}
\zeta =\frac{n}{2}\, f(z)=- i\, \frac{n}{2}\, \log(z/a)\,,
\end{equation}
(we omit the explicit reference to the dependence of $\zeta $ from
$z$, $a$ and $n$ in the notation). Matching condition can be
rewritten as
\begin{equation}\label{firstApproxToE}
E(z)=\left[ I+\mathcal O \left( \frac{1}{n}\right)\right]\, N(z)
\lambda (\beta ;z)^{\sigma _3}\Psi \left( \beta ;  \zeta
\right)^{-1}\,, \quad z\in \CC\,.
\end{equation}
The key idea is to replace $\Psi(\beta ;  \zeta)$ by its leading
asymptotic term as $\zeta \to \infty$. The complete expansion at
infinity of the entries of $\Psi(\beta ;\zeta )$ is well known (see
e.g. \cite[Chapter 9]{abramowitz/stegun:1972}), so we can insert it
in formulas \eqref{RHPPSIx0solution1}--\eqref{RHPPSIx0solution8}.
This computation has been carried out in \cite[Section
4.3]{Vanlessen03}; we can formulate the result therein by defining
the matrix-valued function $\mathcal G$:
\begin{align}
\mathcal G   (\beta ; \zeta )& \isdef
    e^{\frac{\pi i}{4}\sigma_3} e^{-i\zeta\sigma_3} e^{-\frac{1}{2}\beta  \pi
    i\sigma_3}, \quad \text{if $\zeta$ is in the first quadrant,} \label{G0FirstQ} \\
    & \isdef   e^{\frac{\pi i}{4}\sigma_3} e^{-i\zeta\sigma_3} e^{\frac{1}{2}\beta  \pi
    i\sigma_3}, \quad \text{if $\zeta$ is in the second quadrant,} \label{G0SecondQ}  \\
    & \isdef e^{\frac{\pi i}{4}\sigma_3} e^{-i\zeta\sigma_3} e^{\frac{1}{2}\beta  \pi i \sigma_3}
    \begin{pmatrix}
        0 & -1 \\
        1 & 0
    \end{pmatrix}, \quad \text{if $\zeta$ is in the third
    quadrant,} \label{GThirdQ}
    \\ & \isdef  e^{\frac{\pi i}{4}\sigma_3} e^{-i\zeta\sigma_3} e^{-\frac{1}{2}\beta  \pi i \sigma_3}
    \begin{pmatrix}
        0 & -1 \\
        1 & 0
    \end{pmatrix},  \quad \text{if $\zeta$ is in the fourth
    quadrant.} \label{GFourthQ}
\end{align}
It is easy to check that $\mathcal G $ is holomorphic in each
quadrant and matches the jumps of $\Psi  $ on $\Sigma_2$,
$\Sigma_4$, $\Sigma_6$ and $\Sigma_8$ (solid lines in Figure
\ref{fig:local2}), given in \eqref{JumpsForPsi1} and
\eqref{JumpsForPsi4}.
\begin{lemma}[\cite{Vanlessen03}] \label{lemma:completeExpansion} For function
$\Psi(\beta ;\zeta)$ defined in
\eqref{RHPPSIx0solution1}--\eqref{RHPPSIx0solution8},
\begin{equation}\label{AsymptoticsPsiComplete}
\begin{split}
    \Psi(\beta ;\zeta)&=  \frac{1}{\sqrt 2}\,
            \begin{pmatrix}
                1 & -i \\
                -i & 1
            \end{pmatrix}
            \left[
            I+\sum_{k=1}^\infty \frac{i^k}{2^{k+1} \zeta^k}
            \begin{pmatrix}
                (-1)^k s_{\beta ,k} & -i t_{\beta ,k} \\[1ex]
                i (-1)^k t_{\beta ,k} & s_{\beta ,k}
            \end{pmatrix}\right] \,
            \mathcal G(\beta ; \zeta ) ,
\end{split}
\end{equation}
as $\zeta\to\infty$, uniformly for $\zeta$ in each quadrant. Here,
the constants $s_{\beta ,k}$ and $t_{\beta ,k}$ are given by
\begin{equation}\label{sgammatgamma}
        s_{\beta ,k}=\left(\beta +\frac{1}{2},k\right )+\left(\beta -\frac{1}{2},k\right ),\qquad
        t_{\beta ,k}=\left(\beta +\frac{1}{2},k\right )-\left(\beta -\frac{1}{2},k\right ),
\end{equation}
where
\[
    (\nu,k)=\frac{(4\nu^2-1)(4\nu^2-9)\ldots(4\nu^2-(2k-1)^2)}{2^{2k}k!}.
\]
In particular,
\begin{equation}\label{AsymptoticsPsi}
\begin{split}
    \Psi   (\beta ;\zeta)
      &= \frac{1}{\sqrt{2}}  \, \begin{pmatrix}
         1   &  -i   \\[2ex]
        - i   &   1
    \end{pmatrix}\, \left[I -\frac{i \beta   }{2 \zeta} \, \begin{pmatrix}
           \beta      &  i       \\[2ex]
         i  & - \beta
    \end{pmatrix}+\mathcal O\left(\frac{1}{\zeta^2}\right)\right]
      \mathcal G(\beta ; \zeta ) , \quad \zeta \to \infty\,,
\end{split}
\end{equation}
uniformly for $\zeta$ in each sector.
\end{lemma}

Taking into account \eqref{firstApproxToE} and
\eqref{AsymptoticsPsi}, it is reasonable to set in \eqref{formOfP}
\begin{equation}\label{defEfinal}
E(z )   \isdef N(z) \lambda( \beta ; z)^{ \sigma_3} \,   \left[
\frac{1}{\sqrt 2}\,
    \begin{pmatrix}
        1 & -i \\
        -i & 1
    \end{pmatrix}\, \mathcal G\left(\beta ; \zeta   \right)\right]^{-1}\,.
\end{equation}

\begin{proposition}
\label{local_lemma1} Matrix valued function $E$ defined in
\eqref{defEfinal} is holomorphic in $\BB  $ and has there the
following representation:
\begin{equation}\label{representation for E}
    E(z)=\left( \frac{\widehat \FF  (W;z)}{\tau ^2}\, i a^n \right)^{\sigma _3/2}\frac{1}{\sqrt 2}\,
    \begin{pmatrix}
        i & 1 \\
        -1 & -i
    \end{pmatrix}\,,
\end{equation}
where  we take the principal branch of the square root.
\end{proposition}
\begin{proof}
The verification of formula \eqref{representation for E} is
straightforward (and analyticity of $E$ is a direct consequence of
it). Under transformation \eqref{def_zeta_local_analysis} we have
the following correspondence:
\begin{equation*}%\label{correspondence}
    \begin{split}
z \in  (\widehat  \Omega_{\rm i}  \cup \widehat \Omega_0   ) \cap
 \BB^+ \quad & \Leftrightarrow \quad \zeta \text{  is in the first quadrant}, \\
 z \in  (\widehat \Omega_{\rm i}  \cup \widehat \Omega_0   ) \cap
 \BB^- \quad & \Leftrightarrow \quad \zeta \text{  is in the second
 quadrant},\\
 z \in  (\widehat \Omega_{\rm e}  \cup \widehat \Omega_\infty   ) \cap
 \BB^- \quad & \Leftrightarrow \quad \zeta \text{  is in the third
 quadrant},\\
    z \in  (\widehat \Omega_{\rm e}  \cup \widehat \Omega_\infty   )
    \cap \BB^+ \quad & \Leftrightarrow \quad \zeta \text{  is in the fourth
quadrant}.
    \end{split}
\end{equation*}
Assume for instance $z \in  (\widehat  \Omega_{\rm i}  \cup \widehat
\Omega_0   ) \cap  \BB^+$, in which case, according to
\eqref{equ:defN} and \eqref{lambda},
$$
N(z) =\left( \dfrac{D_{\rm e}(W; z)}{\tau} \right)^{\sigma _3}\,,
\quad \lambda( \beta ; z)=  e^{-\pi i \beta } W^{1/2}(z) z^{-n/2}\,,
$$
and by \eqref{G0FirstQ},
$$
\mathcal G   (\beta ; \zeta )=
    e^{\frac{\pi i}{4}\sigma_3} e^{-i\zeta\sigma_3} e^{-\frac{1}{2}\beta  \pi
    i\sigma_3}\,.
$$
Gathering these elements in \eqref{defEfinal} we arrive at
\eqref{representation for E}. The analysis in the rest of the
quadrants is similar.

\end{proof}
\begin{corollary}
Matrix $P(z)=P(a,\beta ;z)$,
\begin{equation}\label{def_P_for_localA}
    P(a,\beta ;z)\isdef \left( \frac{\widehat \FF  (W;z)}{\tau ^2}\,
    i a^n \right)^{\sigma _3/2}\frac{1}{\sqrt 2}\,
    \begin{pmatrix}
        i & 1 \\
        -1 & -i
    \end{pmatrix}\, \Psi \left( \beta ;  \zeta  \right) \,
\lambda (  \beta ; z)^{ - \sigma _3}\,,
\end{equation}
with $\zeta $ given by \eqref{def_zeta_local_analysis}, solves the
Riemann-Hilbert problem (P1)--(P2) defined at the beginning of this
subsection.
\end{corollary}

Thus, it remains to check the matching condition (P3).

\begin{proposition}
\label{local_lemma2} Let $P(z)=P(a,\beta ;z)$ be given by
\eqref{def_P_for_localA}. Then for $z \in \CC $, function $P(z)
N^{-1}(z)$ has the following asymptotic expansion:
\begin{equation}\label{matching_condComplete}
\begin{split}
 P(z) & N^{-1}(z)= I  \\ & +\sum_{k=1}^\infty \frac{i^k}{2^{k+1} \zeta^k}
            \begin{pmatrix}
                s_{\beta ,k} &   (-1)^k  \tau ^{-2} \,  a^n \widehat \FF  (W;z) \,  t_{\beta ,k} \\[1ex]
                \tau ^{2} \,  a^{-n} \left(\widehat \FF  (W;z) \right)^{-1}  t_{\beta ,k} & (-1)^k s_{\beta ,k}
            \end{pmatrix}\,,
\end{split}
\end{equation}
where constants $s_{\beta ,k}$ and $t_{\beta ,k}$ were defined in
\eqref{sgammatgamma}. In particular,
\begin{equation}\label{matching_cond}
\begin{split}
 P(z) N^{-1}(z)=I + \frac{i \beta   }{2 \zeta }   \,
\begin{pmatrix}
             \beta      &      - \tau ^{-2}    a^n \widehat \FF  (W;z)  \\[2ex]
          \tau ^{2}      a^{-n}\left(\widehat \FF (W;z) \right)^{-1}   &  - \beta
    \end{pmatrix}   +  \mathcal O\left(\frac{1}{n^2}
    \right)\,,
\end{split}
\end{equation}
so that $P(z) N^{-1}(z)=I+\mathcal O (n^{-1})$ for $z\in \CC$.
\end{proposition}
\begin{proof}
By \eqref{formOfP},
\[
\begin{split}
P(z) N^{-1}(z) = & E(z)\, \Psi \left( \beta ; \zeta \right) \,
\left[ N(z) \lambda( \beta ; z)^{ \sigma_3} \right]^{-1} \,.
\end{split}
\]
By definition of $E$ in \eqref{defEfinal} we have
$$
N(z) \lambda( \beta ; z)^{ \sigma_3}= E(z ) \left[ \frac{1}{\sqrt
2}\,
    \begin{pmatrix}
        1 & -i \\
        -i & 1
    \end{pmatrix}\, \mathcal G\left(\beta ; \zeta   \right)\right]
    \,.
$$
Thus,
\begin{equation}\label{intermediatePN}
P(z) N^{-1}(z) = E(z)\, \Psi \left( \beta ; \zeta \right) \, \left[
\frac{1}{\sqrt 2}\,
    \begin{pmatrix}
        1 & -i \\
        -i & 1
    \end{pmatrix}\, \mathcal G\left(\beta ; \zeta
    \right)\right]^{-1} E(z)^{-1}.
\end{equation}
By \eqref{AsymptoticsPsi},
\[
\begin{split}
\Psi \left( \beta ; \zeta \right) \, \left[ \frac{1}{\sqrt 2}\,
    \begin{pmatrix}
        1 & -i \\
        -i & 1
    \end{pmatrix}\, \mathcal G\left(\beta ; \zeta
    \right)\right]^{-1} \\ =  \frac{1}{\sqrt 2}\,
    \begin{pmatrix}
        1 & -i \\
        -i & 1
    \end{pmatrix}\,\left[I +V(\zeta )\right] \mathcal G\left(\beta ; \zeta
    \right) \left[ \frac{1}{\sqrt 2}\,
    \begin{pmatrix}
        1 & -i \\
        -i & 1
    \end{pmatrix}\, \mathcal G\left(\beta ; \zeta
    \right)\right]^{-1} \\= \frac{1}{\sqrt 2}\,
    \begin{pmatrix}
        1 & -i \\
        -i & 1
    \end{pmatrix}\,\left[I +V(\zeta )\right]   \left[ \frac{1}{\sqrt 2}\,
    \begin{pmatrix}
        1 & -i \\
        -i & 1
    \end{pmatrix} \right]^{-1},
\end{split}
\]
where $V$ is given by the asymptotic expansion
$$
V(\zeta )\isdef \sum_{k=1}^\infty \frac{i^k}{2^{k+1} \zeta^k}
            \begin{pmatrix}
                (-1)^k s_{\beta ,k} & -i t_{\beta ,k} \\[1ex]
                i (-1)^k t_{\beta ,k} & s_{\beta ,k}
            \end{pmatrix}\,,
$$
and constants $s_{\beta ,k}$ and $t_{\beta ,k}$ were defined
 in \eqref{sgammatgamma}. Using it in \eqref{intermediatePN} and taking into account
\eqref{representation for E} we get
\[
\begin{split}
P(z) N^{-1}(z) =  \left( E(z) \, \frac{1}{\sqrt 2}\,
    \begin{pmatrix}
        1 & -i \\
        -i & 1
    \end{pmatrix} \right)\,\left[I +V(\zeta )\right]   \left( E(z) \, \frac{1}{\sqrt 2}\,
    \begin{pmatrix}
        1 & -i \\
        -i & 1
    \end{pmatrix} \right)^{-1} \\ =
\left( \left( \frac{\widehat \FF (W;z)}{\tau ^2}\, i a^n
\right)^{\sigma _3/2}
    \begin{pmatrix}
        0 & 1 \\
        -1 & 0
    \end{pmatrix} \right)\,\left[I +V(\zeta )\right]   \left( \left( \frac{\widehat \FF (W;z)}{\tau ^2}\, i a^n
\right)^{\sigma _3/2}
    \begin{pmatrix}
        0 & 1 \\
        -1 & 0
    \end{pmatrix} \right)^{-1},
\end{split}
\]
and formula \eqref{matching_condComplete} follows. Furthermore,
since $s_{\beta ,1}=2\beta ^2$, and $t_{\beta ,1}=2 \beta $, we
obtain \eqref{matching_cond}.

\end{proof}

\subsection{Final transformation}

With the notation introduced in \eqref{notationBC} and with
$P(a,\beta ;z)$ defined by \eqref{def_P_for_localA} let us take
\begin{equation*}%\label{def_P_global}
P(z)\isdef P (a_k, \beta _k; z )\quad  \text{for } z \in \BB_k
\setminus (\T  \cup \gamma_{\rm e}\cup \gamma_{\rm i}), \quad k =1,
\dots, m\,,
\end{equation*}
and put
\begin{equation}\label{matrixS}
S(z) \isdef \begin{cases} U(z)N^{-1}(z), & \text{for } z \in
\C\setminus (B \cup  \T  \cup \gamma_{\rm e}\cup \gamma_{\rm i}), \\
U(z)P^{-1}(z), & \text{for } z \in B  \setminus (\T  \cup
\gamma_{\rm e}\cup \gamma_{\rm i}).
\end{cases}
\end{equation}

%%%%%%%%%%%%%%%%%%%%%%%%%%%%%%%%%%%%%%%%%%%%%%%%%%%%%%%%%%
\begin{figure}[htb]
\centering \begin{overpic}[scale=0.7]{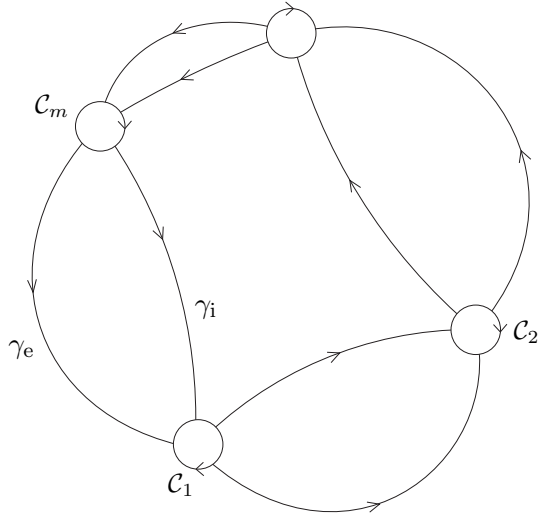}
         \put(35,45){$\gamma_{\rm i}$}
         \put(7,39){$\gamma_{\rm e}$}
          \put(11,75){\small $\CC_m$}
          \put(31,18){\small $\CC_1$}
         \put(83,41){\small $\CC_2$}
\end{overpic}
\caption{Jumps of $S$.}\label{fig:lensesFinal}
\end{figure}
%%%%%%%%%%%%%%%%%%%%%%%%%%%%%%%%%%%%

It is easy to show that this transformation is well defined, since
the inverses exist. Matrix $S$ is holomorphic in the whole plane cut
along $  \gamma \cup C$, where
$$
\gamma  \isdef (\gamma_{\rm e}\cup \gamma_{\rm i}) \setminus B \quad
\text{and} \quad C \isdef \cup_{k=1}^m \CC_k
$$
(see Fig. \ref{fig:lensesFinal}), $S(z) \to I$ as $z\to \infty$, and
if we orient all $\CC_k$'s clockwise, $S_+(t)=S_-(t) J_S$, with
$$
J_S(t)=\begin{cases} P (z) N^{-1}(z), & \text{if } z\in C,  \\
\begin{pmatrix}
1 & 0 \\ \tau^2 /(z^n \FF(W;z)) & 1
\end{pmatrix}, & \text{if } z\in \gamma_{\rm e}\setminus B, \\
\begin{pmatrix}
1 &  - z^n \FF(W; z)/\tau^2  \\ 0 & 1
\end{pmatrix}, & \text{if } z\in \gamma_{\rm i}\setminus B. \\
\end{cases}
$$
It is clear that the off-diagonal terms of $J_S$ on $\gamma_{\rm
i}\setminus B$ and $\gamma_{\rm e}\setminus B$ decay exponentially
fast. On the other hand, by \eqref{matching_cond},
$J_S(z)=I+\mathcal O(1/n)$ for $z\in C$. So the conclusion is that
the jump matrix $J_S=I + O(1/n)$ uniformly for $z \in  \gamma \cup C
$. In fact, \eqref{matching_condComplete} gives us the complete
asymptotic expansion of $J_S$ in negative powers of $n$. Then
arguments such as in
\cite{MR2001f:42037,MR2001g:42050,MR2000g:47048} (see e.g.
\cite[Section 7.2]{MR2001g:42050}) allow to show that $S$ itself has
an asymptotic expansion in negative powers of $n$. The main
observation is that if we define
\begin{equation}\label{def:operatorM}
\mathbf K(f)(z) \isdef \frac{1}{2\pi i}\, \int \frac{ f(t)
(J_S(t)-I)\, dt}{t-z}\,, \qquad \mathbf M(f)(t)=(\mathbf
K(f))_-(t)\,,
\end{equation}
where we integrate along contours $  \gamma \cup C $ with the
orientation shown in Fig.\ \ref{fig:lensesFinal}, then  $ \mathbf M$
defines a bounded linear operator acting in $L^2(\gamma \cup C)$
(with respect to the Lebesgue measure), with the operator norm
$$
\|\mathbf M\|_{L^2\to L^2}\leq \const \|J_S-I\|_{L^2}\,,
$$
where ``$\const$'' depends on $ \gamma \cup C $ only (see e.g.\
\cite[Chapter 7]{MR2000g:47048} for details). In particular, for all
sufficiently large $n$, operator $(1-\mathbf M)$ is invertible, and
we have the following
\begin{proposition}
\label{prop:deift} Let $n$ be such that $(1-\mathbf M)$ is
invertible, and denote
$$
\mu =(1-\mathbf M)^{-1} I\,,
$$
where $I$ is the $2\times 2$ identity matrix. Then matrix $S$ in
\eqref{matrixS} can be expressed as
\begin{equation}\label{equ:exprForSDeift}
S = I +\mathbf K \left(\mu   (J_S(t)-I) \right) \,,
\end{equation}
where we integrate along contours $  \gamma \cup C $ with the
orientation shown in Fig.\ \ref{fig:lensesFinal}.
\end{proposition}

Since for all $n$ sufficiently large, $\|\mathbf M\|_{L^2\to
L^2}<1$, the auxiliary function $  \mu$ can be computed in terms of
the Neumann series,
$$
\mu = \sum_{k=0}^\infty \mathbf M^k (I)\,,
$$
so that
\begin{equation}\label{equ:exprForSDeiftAlt}
S=\sum_{j=0}^\infty S^{(j)}\,, \quad \text{with} \quad S^{(0)}=I\,,
\quad S^{(j+1)}=\mathbf K\left( S^{(j)}_- \right) \,, \quad j\in
\N\,.
\end{equation}
In other words, $S$ can be recovered from $J_S$ by nested contour
integration.

Second observation is that, taking into account the exponential
decay of the off-diagonal terms of $J_S$ on $\gamma_{\rm e}\setminus
B$ and $ \gamma_{\rm i}\setminus B$, we can restrict the integration
in \eqref{def:operatorM} to contours $\CC$, replacing then $J_S$ by
$P N^{-1}$; in this way, only exponentially small terms are
neglected. Pluggin the asymptotic expansion of $P N^{-1}$, given by
\eqref{matching_condComplete}, into equation
\eqref{equ:exprForSDeiftAlt} allows to find successively the terms
$\mathfrak s_k$. In particular, for terms $S^{(j)}$ in
\eqref{equ:exprForSDeiftAlt} we have
$$
S^{(j)}(z)=\mathcal O\left( \frac{1}{n^j}\right)
$$
locally uniformly in $\C$. Furthermore,
\begin{equation}\label{negativeS}
    S(z)=I+\sum_{k=1}^\infty \frac{\mathfrak s_k(z)}{n^k}\,,
\end{equation}
uniformly in $z$, where $\mathfrak s_k$'s are piece-wise analytic
functions in $\C \setminus \CC$. Let us determine explicitly the
first nontrivial term $\mathfrak s_1$. By \eqref{def:operatorM} and
\eqref{equ:exprForSDeiftAlt},
$$
S^{(1)}(z)= \mathbf K\left( I \right)=- \sum_{k=1}^m \frac{1}{2\pi
i}\, \oint_{\CC_k} \frac{ (P(t) N^{-1}(t)-I)\, dt}{t-z}
+\text{exponentially small terms,}
$$
where integrals are taken counterclockwise. Hence, we have
\begin{corollary}
Matrix $S$ satisfies
\begin{equation*} %\label{asympS}
    S(z) = I - \sum_{k=1}^m \frac{1}{2\pi i}\, \oint_{\CC_k} \frac{ (P(t) N^{-1}(t)-I)\,
dt}{t-z} +   \mathcal O\left(\frac{1}{n^2}\right)
    \end{equation*}
locally uniformly for $z \in \C\setminus ( \gamma \cup C)$, where
integrals along $\CC_k$'s are taken counterclockwise.
\end{corollary}

Define for $z\in \BB_k \cup \CC_k$ ($k=1, \dots, m$)
\begin{equation*}%\label{def_Hat_S}
    \mathcal  F_k(z)\isdef - \frac{  \beta_k   }{n \log (z/a_k)}   \,
\begin{pmatrix}
             \beta_k      &      - \tau ^{-2}    a_k^n\, \widehat
             \FF_k
    (W;z)       \\[2ex]
          \tau ^{2}       a_k^{-n}\, \widehat \FF_k^{-1}
    (W;z)   &  - \beta_k
    \end{pmatrix}\,,
\end{equation*}
so that by \eqref{matching_cond},
$$
P(z)N^{-1}(z)-I=\mathcal F_k(z)+\mathcal
O\left(\frac{1}{n^2}\right), \quad z \in \CC_k\,.
$$
Since $ \mathcal  F_k(z)$ is a meromorphic function in $\BB_k \cup
\CC_k$ with a simple pole at $z=a_k$, by the residue theorem, for $z
\notin \BB_k \cup \CC_k$,
$$
\frac{1}{2\pi i}\, \oint_{\CC_k } \frac{ \mathcal F_k(t)\,
dt}{t-z}=\res_{t=a_k}   \frac{ \mathcal F_k(t)\, dt}{t-z}=F_k(z)\,,
$$
where
\begin{equation*}%\label{def_residues}
F_k(z)\isdef -\dfrac{a_k \beta _k  }{n (a_k-z)} \,
\begin{pmatrix}
             \beta_k      &      - \tau ^{-2} a_k^n \vartheta_k    \\[2ex]
         \tau ^{ 2}  a_k^{-n} \vartheta_k^{-1}  &  - \beta_k
    \end{pmatrix}
   \,.
\end{equation*}
Consequently,
$$
\frac{1}{2\pi i}\, \oint_{\CC_k} \frac{ (P(t) N^{-1}(t)-I)\,
dt}{t-z} =\begin{cases}    F_k(z) +\mathcal
O(1/n^2), & \text{if } z \in \C \setminus \overline{\BB_k}\,, \\
     (\mathcal F_k + F_k)(z) +\mathcal O(1/n^2), & \text{if } z \in  \BB_k\,.
\end{cases}
$$
In particular, in \eqref{negativeS},
\begin{equation}\label{residues_Final}
\mathfrak s_1(z)=    \sum_{k=1}^m  \dfrac{a_k \beta _k }{
 a_k-z } \,
\begin{pmatrix}
             \beta_k      &      - \tau ^{-2} a_k^n \vartheta_k    \\[2ex]
         \tau ^{ 2}  a_k^{-n} \vartheta_k^{-1}  &  - \beta_k
    \end{pmatrix}   \quad \text{for } z \in \C\setminus
(\gamma\cup B)\,.
\end{equation}

Now we are ready for the asymptotic analysis of the original matrix
$Y$ (and in particular, of its entries $(1,1)$ and $(2,1)$), that we
perform in the next section.

\section{Asymptotic analysis} \label{sec:asymptotics}

Unraveling our transformations we have
\begin{equation}\label{equ:exprForYlast3}
    Y(z)=\begin{cases}
     S(z)N(z)K^{-1}(z)\, H^{-1}(z), & \text{if } z \in \C \setminus B , \\  S(z)
P(z)K^{-1}(z)\, H^{-1}(z), & \text{if } z \in  B.
    \end{cases}
\end{equation}

%%%%%%%%%%%%%%%%%%%%%%%%%%%%%%%%%%%%%%%%%%%%%%%%%%%%%%%%%%
\begin{figure}[htb]
\centering \begin{overpic}[scale=0.7]{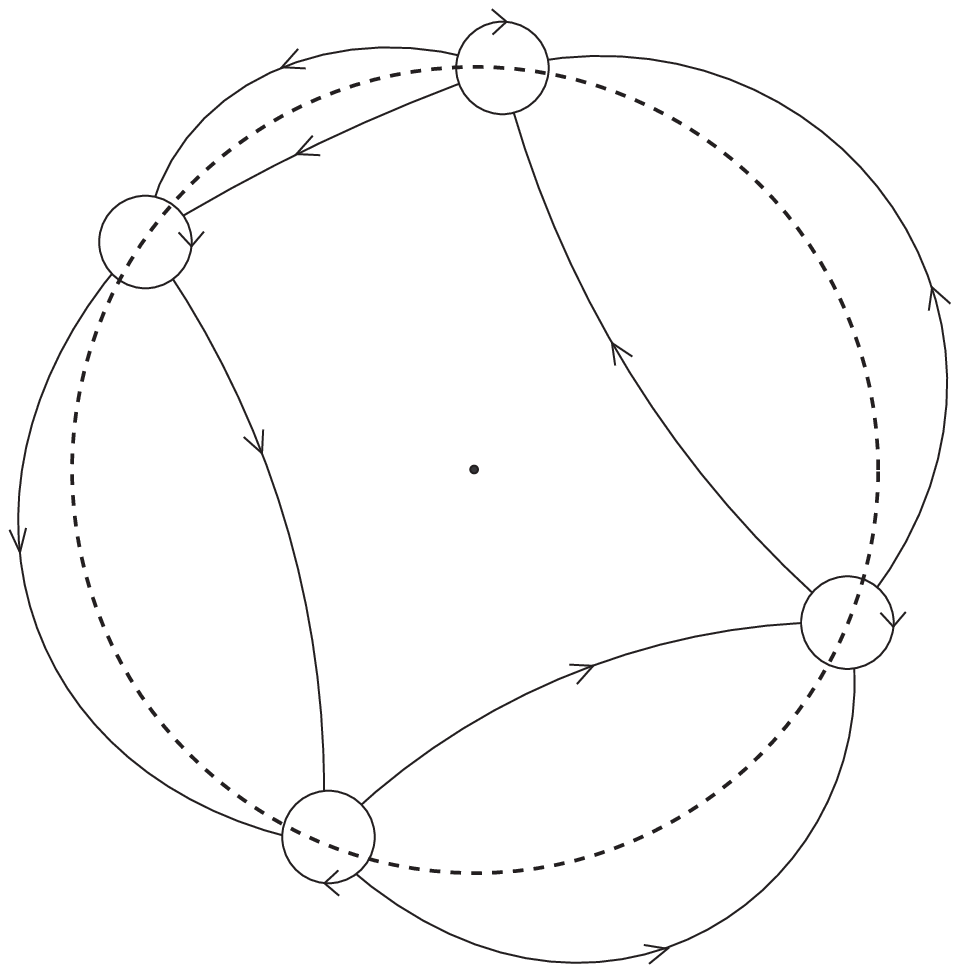}
      \put(55,27){\small $\T $}
         \put(35,45){$\gamma_{\rm i}$}
         \put(7,39){$\gamma_{\rm e}$}
          \put(40,65){$\Omega_0$}
          \put(65,69){$\Omega_{\rm i}$}
          \put(74,75){$\Omega_{\rm e}$}
          \put(87,82){$\Omega_\infty$}
          \put(11,75){\small $c_m$}
          \put(32,19){\small $c_1$}
         \put(82,41){\small $c_2$}
\end{overpic}
\caption{Domains for the asymptotic
analysis.}\label{fig:lensesFinalAsymptotics}
\end{figure}
%%%%%%%%%%%%%%%%%%%%%%%%%%%%%%%%%%%%

\begin{proof}[Proof of Theorem \ref{thm1} and Corollary \ref{cor:zeros}]
Using \eqref{negativeS} in \eqref{equ:exprForYlast3} we readily
obtain \eqref{completeExpansionofPhi}. In order to prove the rest of
the statements we must analyze the consequences of
\eqref{equ:exprForYlast3} in each domain (see Fig.\
\ref{fig:lensesFinalAsymptotics}).

In $\Omega_0\setminus B$  we have (cf.\ \eqref{defH}, \eqref{KforT1}
and \eqref{equ:defN})
$$
N(z)=    \begin{pmatrix} 0 & D_{\rm i}(W; z)/\tau   \\ -\tau /D_{\rm
i}(W; z) & 0
\end{pmatrix}, \quad K(z)=H(z)=I\,.
$$
Hence, $  Y(z)=S(z) N(z)$, so that
$$
Y_{11}(z)=-\frac{\tau }{D_{\rm i}(W; z)}\, S_{12}(z)\,, \quad
Y_{21}(z)= -\frac{\tau }{D_{\rm i}(W; z)}\, S_{22}(z)\,.
$$
Taking into account \eqref{residues_Final} and recalling that
$\Phi_n =Y_{11} $ we obtain
\begin{equation*}%\label{asymptForPolynsInSigma0}
\Phi_n(z)=  \frac{1}{\tau \, D_{\rm i}(W; z)}\, \frac{1}{n}\, \left(
\sum_{k=1}^m \dfrac{\beta _k \vartheta_k a_k^{n+1}  }{ a_k-z }   +
O\left(\frac{1}{n }\right) \right)\,, \quad z \in \Omega_0\setminus
B\,,
\end{equation*}
which yields formula (i) in Theorem \ref{thm1}.

Observe that the leading term in the right hand side is a rational
function with at most $m-1$ zeros. In consequence, for all
sufficiently large $n$ each $\Phi_n(z)$ can have at most $m-1$ zeros
on a compact subset of the unit disk (Corollary \ref{cor:zeros}).

Let us analyze the asymptotic behavior of the polynomials close to
the inner boundary of the unit circle, but still away from the
singular points $\mathcal A$. For $z \in \Omega_{\rm i}\setminus B$
we have
\begin{equation}\label{NKH_case1}
N(z)=    \begin{pmatrix} 0 & D_{\rm i}(W; z)/\tau   \\ -\tau /D_{\rm
i}(W; z) & 0
\end{pmatrix}, \quad K^{-1}(z)=\begin{pmatrix} 1 & 0  \\ z^{n}/  W(z)  & 1
\end{pmatrix}\,,\quad  H(z)=I \,.
\end{equation}
Hence,
\[
 \begin{split}
  Y(z)&=S(z)  \,   \begin{pmatrix}  z^{n}D_{\rm
i}(W; z) /(\tau W(z)) & D_{\rm i}(W; z)/\tau
\\ -\tau /D_{\rm i}(W;z)     & 0
\end{pmatrix}=S(z)  \,   \begin{pmatrix}  z^{n}D_{\rm
e}(W; z) / \tau   & D_{\rm i}(W; z)/\tau
\\ -\tau /D_{\rm i}(W;z)     & 0
\end{pmatrix}  . \end{split}
\]
Analyzing $Y_{11}$ we obtain that
\begin{equation}\label{asympExprinOmegaPlus}
 \begin{split}
\Phi_n(z) = &  \frac{z^n D_{\rm e}(W;z)}{\tau }\, \left(1+
\frac{1}{n}\, \sum_{k=1}^m \frac{a_k \beta _k^2  }{ a_k-z }
+\mathcal O \left(\frac{1}{n^2}\right)\right) \\ & + \frac{1 }{
\tau\, D_{\rm i}(W;z)}\, \left(\frac{1}{n}\, \sum_{k=1}^m \frac{
\beta _k \vartheta_k a_k^{n+1}}{ a_k-z }  +\mathcal O
\left(\frac{1}{n^2}\right) \right)\,,
 \end{split}
\end{equation}
which coincides with the formula of statement (ii) of Theorem
\ref{thm1}.

\medskip

The exterior asymptotic is analyzed likewise. For $z\in \Omega_{\rm
e}\setminus B $  we have (cf.\ \eqref{defH}, \eqref{KforT1} and
\eqref{equ:defN})
$$
N(z)=  (\left( D_{\rm e} (W;z)/\tau  \right)^{ \sigma_3} , \quad
K^{-1}(z)=\begin{pmatrix} 1 & 0
\\ -1/(z^{n } W(z)) & 1
\end{pmatrix}, \quad H^{-1}(z)=z^{n \sigma _3}\,.
$$
Unraveling the transformation we conclude that formula
\eqref{asympExprinOmegaPlus} is valid locally uniformly also in
$\Omega_{\rm e}\setminus B $, that is, holds in a neighborhood of
$\T $ away from the singular set $\AA$.

For $z \in \Omega_\infty \setminus B$  we have
$$
N(z)=  (\left( D_{\rm e} (W;z)/\tau  \right)^{ \sigma_3} , \quad K
=I, \quad H^{-1}(z)=z^{n \sigma _3}\,.
$$
Hence, $ Y(z) =S(z) \, (\left( z^n D_{\rm e} (W;z)/\tau  \right)^{
\sigma_3} $, and statement (iii) of Theorem \ref{thm1} is a direct
consequence of formula \eqref{residues_Final}.
\end{proof}

\begin{proof}[Proof of Theorem \ref{thm:local_behavior}]
Now we turn to the asymptotic analysis for $\Phi_n$'s in a
neighborhood of a singularity $a_k$, where we have to use the
expression $$Y(z)= S(z) P(z)K^{-1}(z)\, H^{-1}(z)$$ with $P(z)= P
(a_k, \beta _k; z )$ for $z \in \BB_k$.

Assume that $z  \in \Omega_{\rm i}\cap \BB_k $ and $\arg(z) >
\arg(a_k)$ in such a way that $\zeta= -i \,\frac{n}{2}\,
\log\left(z/a_k \right)$ belongs to sector I in Figure
\ref{fig:local2}. Then   $\Psi$ is given by
\eqref{RHPPSIx0solution1}, which along with \eqref{lambda} replaced
in \eqref{def_P_for_localA} yields
$$
P(z)=\frac{\sqrt\pi}{2 \sqrt 2} \, \zeta^{1/2} \left( \frac{\widehat
S_k(W;z)}{\tau ^2}\, i a_k^n \right)^{\sigma _3/2} M(\zeta ) \left(
e^{-\pi i /4} W^{-1/2}(z) z^{n/2} \right)^{\sigma _3},
$$
where
\begin{equation*}%\label{def_M}
M(\zeta  ) \isdef
    \begin{pmatrix}
      i  H_{\beta_k  +\frac{1}{2}}^{(2)}+H_{\beta_k  -\frac{1}{2}}^{(2)}  &
         H_{\beta_k  +\frac{1}{2}}^{(1)}   -i H_{\beta_k  -\frac{1}{2}}^{(1)}  \\[2ex]
       -H_{\beta_k  +\frac{1}{2}}^{(2)} - i H_{\beta_k  -\frac{1}{2}}^{(2)}  &
       i  H_{\beta_k  +\frac{1}{2}}^{(1)}   -  H_{\beta_k
       -\frac{1}{2}}^{(1)}
    \end{pmatrix}(\zeta)
    \,.
\end{equation*}
Thus,
\begin{align*}
P_{11}(z) &= \frac{\sqrt{\pi}}{2 \sqrt{2}} \,   \zeta^{1/2}
 \left( \dfrac{\widehat
S_k(W;z)}{\tau ^2}\,  a_k^n \frac{z^n}{W(z)} \right)^{1/2} \, M_{11}(\zeta ) \,, \\
P_{12}(z) &=   \frac{\sqrt{\pi}}{2 \sqrt{2}} \,   \zeta^{1/2}
  \left( -\dfrac{\widehat
S_k(W;z)}{\tau ^2}\,  a_k^n \frac{W(z)}{z^n} \right)^{1/2}  \,   M_{12}(\zeta ) \,, \\
P_{21}(z) &=   \frac{\sqrt{\pi}}{2 \sqrt{2}} \,   \zeta^{1/2} \,
\left( -\dfrac{\tau ^2}{\widehat
S_k(W;z)}\,  a_k^{-n} \frac{z^n}{W(z)} \right)^{1/2}  \,   M_{21}(\zeta )\,, \\
P_{22}(z) &=   \frac{\sqrt{\pi}}{2 \sqrt{2}} \,   \zeta^{1/2}  \,
\left(  \dfrac{\tau ^2}{\widehat S_k(W;z)}\,  a_k^{-n}
\frac{W(z)}{z^n} \right)^{1/2}  \, M_{22}(\zeta )\,.
\end{align*}
Taking into account that   $K$ and $H$ are as in \eqref{NKH_case1}
we see that
\[
\begin{split}
(P(z) K^{-1}(z))_{11}&= P_{11}(z) +P_{12}(z)\, \frac{z^n}{W(z)}
\\ & =\frac{\sqrt{\pi}}{2 \sqrt{2}} \,   \zeta^{1/2}
  \left( \dfrac{\widehat
S_k(W;z)}{\tau ^2}\,  a_k^n \frac{z^n}{W(z)} \right)^{1/2} \, \left(
M_{11}(\zeta )+i M_{12}(\zeta )\right) \,,
\\
(P(z) K^{-1}(z))_{21} & = P_{21}(z) +P_{22}(z)\, \frac{z^n}{W(z)} \\
&=  \frac{\sqrt{\pi}}{2 \sqrt{2}} \,   \zeta^{1/2}
 \left(   \dfrac{\tau ^2}{\widehat
S_k(W;z)}\,  a_k^{-n} \frac{z^n}{W(z)} \right)^{1/2} \, \left( i
M_{21}(\zeta )+M_{22}(\zeta ) \right)\,.
\end{split}
\]
Thus,
\begin{equation*}
%\label{dominant}
\begin{split}
Y_{11}(z)& =S_{11}(z)(P(z) K^{-1}(z))_{11}+S_{12}(z)(P(z)
K^{-1}(z))_{21} \\ &= (P(z) K^{-1}(z))_{11} \, \left( 1+\mathcal
O\left(\frac{1}{n}\right)\right)\,.
\end{split}
\end{equation*}
By \cite[formulas 9.1.3--9.1.4]{abramowitz/stegun:1972},
$$
M_{11}(\zeta )+i M_{12}(\zeta )=2\left(i J_{\beta_k +\frac{1}{2}}
+J_{\beta_k  -\frac{1}{2}} \right)(\zeta )\,,
$$
so that for $z$ in the domain considered
\begin{equation}\label{asymptoticsNearSingularity1}
\begin{split}   \Phi_n(z) &= \sqrt{\frac{\pi   } { 2}} \,
 \left( \dfrac{\widehat
S_k(W;z)}{\tau ^2}\,  \frac{a_k^n  z^n}{W(z)} \right)^{1/2} \mathcal
H(\beta_k ; \zeta )\, \left( 1+\mathcal
O\left(\frac{1}{n}\right)\right)\,,
\end{split}
\end{equation}
where $\zeta=-i \frac{n}{2}\, \log\left( z/a_k\right)$, and
$\mathcal H(\beta ; \zeta )$ has been defined in
\eqref{def_H_for_local}. Function $\mathcal H(\beta ; \zeta )$ is
holomorphic in the complex plane cut along the positive imaginary
axis $i\R_+$, and by \cite[formula 9.1.35]{abramowitz/stegun:1972},
\begin{equation}\label{jumpsH_across_imaginary}
\mathcal H_+(\beta ; \zeta  )=e^{2\pi i \beta } \mathcal H_-(\beta ;
\zeta )\,, \quad \zeta \in i\R_+\,,
\end{equation}
if the positive imaginary axis is oriented towards the origin.

Taking into account that
$$
\left( \dfrac{\widehat S_k(W;z)}{W(z)}\right)^{1/2}=\begin{cases}
e^{ \pi i \beta _k/2} D_{\rm e}(W; z)\,, & \text{if } z\in \BB_k
\text{ and } \arg(z)>\arg(a_k)\,, \\ e^{- \pi i \beta _k/2} D_{\rm
e}(W; z)\,, & \text{if } z\in \BB_k \text{ and } \arg(z)<
\arg(a_k)\,,
    \end{cases} \quad k=1, \dots, m\,,
$$
we see that \eqref{Local_asymptotics} is valid for $z  \in
\Omega_{\rm i}\cap \BB_k $ and $\arg(z) > \arg(a_k)$.

Since  transformation $z \mapsto \zeta $ maps the ray $a_k\cdot
(0,1)$  onto the positive imaginary axis oriented towards the
origin, from formulas \eqref{boundaryValueD_onCuts} and
\eqref{jumpsH_across_imaginary} it follows that $D_{\rm e}(W;z)
\mathcal H(\beta_k ; \zeta_n )$ is single-valued in a neighborhood
of $z=a_k$. Hence, by uniqueness of the analytic continuation we
must conclude that \eqref{Local_asymptotics} is valid in fact in the
whole $\BB_k$.
\end{proof}

\medskip

\begin{proof}[Proof of Theorem \ref{thm:coefficients}]
For the Verblunsky coefficients we have also
$$
 \overline{\alpha _n}=-\Phi_{n+1}(0)=  - \frac{1}{n}\, \sum_{k=1}^m
\beta _k \vartheta_k a_k^{n+1}  + O\left(\frac{1}{n^2 }  \right) \,,
\quad n \to \infty\,,
$$
and it remains to take into account that $\beta _k\in \R$,
$\vartheta_k \in \T $, and $a_k\in \T $ in order to arrive at
formula \eqref{result_for_Verblunsky}.

Furthermore, observe from \eqref{Y} that the leading coefficients
$\kappa_n$ of the orthonormal polynomials $\varphi_n$ have the
following representation in terms of matrix $Y$:
$$
\kappa_{n-1}^2=-\frac{1}{2\pi}\,Y_{21}(0)\,.
$$
Since in $\Omega_0\setminus B$  we have $  Y(z)=S(z) N(z)$, with
$$
N(z)=    \begin{pmatrix} 0 & D_{\rm i}(W; z)/\tau   \\ -\tau /D_{\rm
i}(W; z) & 0
\end{pmatrix},
$$
then
$$
Y_{21}(z)= -\frac{\tau }{D_{\rm i}(W; z)}\, S_{22}(z)\,,
$$
so that
$$
\kappa_{n-1}^2= \frac{\tau ^2}{2\pi}\,S_{22}(0)\,,
$$
and formula \eqref{result_for_kappa} follows from
\eqref{tauAlternativeBis} and \eqref{residues_Final}.

\end{proof}

\begin{proof}[Proof of Theorem \ref{thm:Fisher}]
One of the main connections of Toeplitz determinants $\mathcal
D_n(W)$ defined in \eqref{def_Toeplitz_dets} with the orthogonal
polynomials on the unit circle is that they can be expressed in
terms of the leading coefficients $\kappa _n$ using the following
formula (see e.g.\ \cite[Theorem 1.5.11]{Simon05b}):
$$
\frac{\mathcal D_n(W)}{\mathcal D_{n-1}(W)}=\frac{1}{\kappa _n^2}\,.
$$
Since $\mathcal D_0(W)=d_0=\oint z^{-n} W(z) |dz|$, we obtain that
\begin{equation*}%\label{toeplitz1}
    \mathcal D_n(W)=d_0\, \prod_{j=1}^n \kappa _j^{-2}\,.
\end{equation*}
Fix $N\in \N$ such that $N >\sum_{k=1}^m \beta _k^2$. By formula
\eqref{result_for_kappa}, for $n>N$,
\[
\begin{split}
 \mathcal D_n(W) &=d_0\, \prod_{j=1}^{N-1} \kappa _j^{-2}
  \prod_{j=N}^n \kappa _j^{-2} \\ &= d_0\, \prod_{j=N}^n \left[  \frac{\tau^2}{2\pi}\,\left(
  1- \frac{1}{j+1}\, \sum_{k=1}^m \beta _k^2+\mathcal O\left(
\frac{1}{j^2}\right) \right) \right]^{-1}
\end{split}
\]
where we denote $\mathcal E_1\isdef d_0\, \prod_{j=1}^{N-1} \kappa
_j^{-2}$. Hence,
\[
\begin{split}
 \mathcal D_n(W)  &= \mathcal E_1 \left(
\frac{2\pi}{\tau^2} \right)^n\, \prod_{j=N}^n \left[  1-
\frac{1}{j+1}\, \sum_{k=1}^m \beta _k^2+\mathcal O\left(
\frac{1}{j^2}\right)   \right]^{-1} \\ &= \mathcal E_1  \left(
\frac{2\pi}{\tau^2} \right)^n\, \prod_{j=N}^n \left[  1-
\frac{1}{j+1}\, \sum_{k=1}^m \beta _k^2   \right]^{-1} \prod_{j=N}^n
\left[1+\frac{ \mathcal O\left(  j^{-2} \right) }{ 1-
\frac{1}{j+1}\, \sum_{k=1}^m \beta _k^2+\mathcal O\left( j^{-2}
\right) }\right]^{-1}  .
\end{split}
\]
This last infinite product is convergent to a constant, which we
denote by $\mathcal E_2$. Hence,
\[
\begin{split}
\mathcal D_n(W) & =\mathcal E_1  \mathcal E_2 \left(
\frac{2\pi}{\tau^2} \right)^n\, \prod_{j=N}^n \left[  1-
\frac{1}{j+1}\, \sum_{k=1}^m \beta _k^2 \right]^{-1} \\ &= \mathcal
E_1  \mathcal E_2 \left( \frac{2\pi}{\tau^2} \right)^n\,
\frac{\Gamma(n+2)\, \Gamma\left(N+1- \sum_{k=1}^m \beta _k^2\right)
}{\Gamma(N+1)\, \Gamma\left(n+2- \sum_{k=1}^m \beta _k^2\right)} \,.
\end{split}
\]
Gathering all the constants in $\varkappa$ and using Stirling
formula for the asymptotics of the Gamma function we obtain that
\[
\begin{split}
\mathcal D_n(W) & =\varkappa\,  \left( \frac{2\pi}{\tau^2}
\right)^n\,  n^{\sum_{k=1}^m \beta _k^2} \, \left( 1+\mathcal O
(1)\right) \,, \quad n \to \infty\,.
\end{split}
\]
It remains to use formula \eqref{tauAlternativeBis} in order to set
the proof of Theorem \ref{thm:Fisher}.
\end{proof}

 \begin{proof}[Proof of Theorem \ref{thm:clock}]
Let  $\varepsilon >0$. We may rewrite formula \emph{(iii)} of
Theorem \ref{thm1} using the notation \eqref{G_n}:
\begin{equation}\label{case_finite_number_dominant_singRewritten}
\frac{  D_{\rm i}(W; z)}{a_1^{n+1} D_{\rm i}(W; 0)}\, \Phi_n(z)=
\frac{D_{\rm i}(W; z)}{a_1^{n+1}}\, \left\{ z^n \mathcal S(W;z)
-\dfrac{1}{n}\, \mathcal R_n(z) \right\} + H_n(z) \,.
\end{equation}
If $z \in \Gamma_n(\varepsilon)$, then there exists a constant
$C=C(\varepsilon)$ such that
\begin{equation}\label{boundForH}
|H_n(z)|\leq \frac{C}{n^2}\,.
\end{equation}
We conclude the proof using standard arguments involving Rouche's
theorem (see e.g.\ proof of Theorem 4 in \cite{math.CA/0502300}).
\end{proof}

\begin{proof}[Proof of Theorem \ref{thm:szabados}]
Let $K\subset \D$ be a compact set. For $n$ large enough and $z\in
K$, we have $|z|^n \leq n^{-2}$, so that by
\eqref{case_finite_number_dominant_singRewritten},
$$
    Z \cap K \subset \bigcap_{k\geq 1} \overline{\bigcup_{n\geq k} \mathcal Z( \mathcal R_n)
    }\,,
$$
with $\mathcal R_n$ defined in \eqref{G_n}. Thus, it is sufficient
to describe all the possible limit points of $\{ \mathcal G_n\} $.
Observe that with the notation introduced in Section
\ref{section:intro}, before the formulation of Theorem
\ref{thm:szabados},
\begin{align*}
 \mathcal R_n(z)& = \sum_{k=1 }^m\frac{\beta _k \vartheta_k  }{z- a_k }\,
 \exp\left(2\pi i (n+1) \theta_k\right) = \sum_{k=1 }^m\frac{\beta _k
 \vartheta_k  }{z- a_k } \, \exp\left(2\pi i  \sum_{j=1}^v r_{kj}\, (n-m+1)
 \theta_j\right)\,.
\end{align*}

Let $v=1$; this means that all $\theta_k\in \Q$, and
$\theta_k=r_{k1}\equiv p_k/q_k \mod \Z$, $k=2, \dots, m$. Then
clearly all possible limits of $\mathcal R_n(z)$ are given by
equation \eqref{description_Szabados}.

Assume $v\geq 2$. By Kronecker's theorem (also known as
Kronecker-Weyl theorem, see e.g.~\cite[Ch.\ III]{Cassels:1957}),
since $\theta_1=1, \theta_2, \dots, \theta_v$ are rationally
independent then (and in fact, if and only if) sequence
$$
\left\{ \left(e^{2\pi i \theta_2 n}, \dots ,  e^{2\pi i \theta_v
n}\right) \right\}_{n\in \N} \subset \T ^{v-1}
$$
is dense and uniformly distributed in the $(v-1)$-dimensional torus
$\T ^{v-1}$. In particular, for any real numbers $X_2, \dots, X _v$
there exists a sub-sequence $\Lambda \subset \N$ such that
\begin{equation}\label{Kronecker1}
\lim_{n\in \Lambda} e^{ 2\pi i   \theta _j n}=e^{ 2\pi i   X _j} \,,
\quad j=2, \dots, v\,.
\end{equation}
In fact, we can say more:
\begin{lemma} \label{lemma:KroneckerAux}
Let $\theta_1=1, \theta_2, \dots, \theta_v$ be rationally
independent, and let $r_{kj}\in \Q$, $j=2,\dots, v$, $k=1, \dots,
m$. Then for any real numbers $X_2, \dots, X _v$ there exists a
sub-sequence $\Lambda \subset \N$ such that
\begin{equation}\label{Kronecker2}
\lim_{n\in \Lambda} e^{ 2\pi i r_{kj}  \theta _j n}=e^{ 2\pi i
r_{kj}  X _j} \,, \quad j=2, \dots, v\,, \quad k=1, \dots, m\,.
\end{equation}
\end{lemma}
Indeed, for $r_{kj}\in \Z$ this statement follows trivially from
\eqref{Kronecker1}. If $1/r_{kj}\in \Z$, then
$$
e^{ 2\pi i r_{kj}  \theta _j n}=\left( e^{ 2\pi i  \theta _j
n}\right)^{\frac{1}{1/r_{kj} }}\,,
$$
and we can specify the single-valued branch of the $r_{kj}^{-1}$-th
root in the neighborhood of $\exp(2\pi i X_j)$ in such a way that
the corresponding limit in \eqref{Kronecker2} holds. Combining these
two observations we obtain Lemma \ref{lemma:KroneckerAux}, which
shows that the set of limit points of $\{ \mathcal R_n\} $ is given
by the left hand side of \eqref{description_Szabados2}. Now the
statement follows for $v\geq 2$.
\end{proof}

\begin{proof}[Proof of Corollary \ref{cor:cluster}]
There is a finite number of numbers $s_k\in [0,q_k) \cap \Z$; for
each possible combination of $s_k$'s, the left hand side in
\eqref{description_Szabados} is a rational function with denominator
of degree $\leq m-1$. Now the statement \emph{(i)} follows.

If $v=2$, then $Z\cap \D$ is a manifold parameterized by a
continuous parameter $X_2\in\R$, which shows that it is a curve. It
is easy to check that its equation is a polynomial in two real
variables of degree $\leq m$. Furthermore, if $v=m=2$, then equation
\eqref{description_Szabados2} is equivalent to
$$
\frac{\beta _1 \vartheta_1 (z-a_2)}{\beta _2
\vartheta_2(a_1-z)}=e^{2\pi i X_2}  \quad \Longleftrightarrow \quad
\left| \frac{\beta _1 \vartheta_1 (z-a_2)}{\beta _2
\vartheta_2(a_1-z)}\right|=1\,,
$$
which reduces to $|\beta _1| \, |z-a_1| =|\beta _2| \, |z-a_2|$.
This proves \emph{(ii)}.

Finally, for $v>2$, the set $Z\cap \D$ is a manifold parameterized
by at least two continuous parameter $X_i\in\R$, showing that
generically it is a two-dimensional domain. Its boundary is again an
algebraic curve or a union of algebraic curves.
\end{proof}

\begin{proof}[Proof of Theorem \ref{thm:closestToA}]
It is a straightforward consequence of the asymptotic formula
\eqref{Local_asymptotics}, that shows that in the neighborhood
$\BB_k$ of $a_k$ the zeros of $\Phi_n$ match (up to a $\mathcal
O(1/n)$ term) those of $\mathcal H(\beta_k ; \zeta_n )$.
\end{proof}

\begin{remark}\label{remark:zero}
We don't know any explicit formula for $h(\beta )$, although it can
be easily computed numerically. In order to find a good initial
value for $h(\beta )$ we can use the continued fraction expansion
\cite[Formula 9.1.73]{abramowitz/stegun:1972}:
$$
 \frac{J_\nu(z)}{J_{\nu-1}(z)} = \dfrac{1}{2\nu /z-\dfrac{1}{2(\nu +1)/z-\dfrac{1}{\ddots}}}\,.
$$
In particular, truncating at the second term and equating to zero we
can take
$$
  h_0(\beta )=\frac{\sqrt{(2\beta +3)(6\beta +1)}+ i (2\beta
+3)}{2}
$$
as a reasonably good approximation for any iterative (say,
Newton-type) zero-finding method of computation of $h(\beta )$. As
an illustration, in the table below we compare the values of
$h_0(\beta )$ and $h(\beta )$ for $\beta =-1/4$ and $\beta =1,
\dots, 5$:
\begin{center}
\begin{tabular}{ |c||c|c|c|c|c|c| } \hline %
$ \beta  $ & $-0.25$ & $1$ & $2$ & $3$ & $4$ & $5$ \\ \hline %
$h_0(\beta )$ & $0.69 i$ & $2.96+2.5 i$ & $4.77+3.5 i$ &
$ 6.54 + 4.5 i$ & $ 8.3 + 5.5 i$ & $ 10.04 + 6.5 i $ \\
\hline $ h(\beta )$ & $0.68 i$ & $ 3.73+ 1.04 i$ & $ 5.08 + 0.87 i
$& $6.34 + 0.79 i$  & $7.56 + 0.74 i$ &
$8.75 + 0.71 i$ \\ \hline %
\end{tabular}
\end{center}
\end{remark}

\section*{Acknowledgement} The research of A.M.F.\ was supported, in part, by
a grant from the Ministry of Education and Science of Spain, project
code MTM2005-08648-C02-01, by Junta de Andaluc\'{\i}a, grants FQM229 and
FQM481, and by ``Research Network on Constructive Complex
Approximation (NeCCA)'', INTAS 03-51-6637.

The research of K.T.-R.M.\ was supported, in part, by the U.S.
National Science Foundation under grants DMS--0200749 and
DMS--0451495.

Both A.M.F.\ and K.T.-R.M.\ acknowledge also a partial support of
the NATO Collaborative Linkage Grant ``Orthogonal Polynomials:
Theory, Applications and Generalizations'', ref. PST.CLG.979738.

The research of E.B.S.\ was supported, in part, by the U.S. National
Science Foundation under grant DMS--0532154.

Finally, we are grateful to Prof.\ Muldoon for his useful comments
included in Remark \ref{remark:localBehavior:new}.

\end{document}